\documentclass[11pt,draftcls,onecolumn]{IEEEtran}

\usepackage{cite}

\usepackage{amsmath,bbm,times,stmaryrd}

 \usepackage[mathscr]{eucal}
 \usepackage{amsmath,amssymb}
\usepackage{tabularx}
\usepackage{stmaryrd}
\usepackage{multirow}
\usepackage{enumerate}

\usepackage[dvipsnames,usenames]{color}
\usepackage{psfrag}
\usepackage[dvips]{graphicx}
\usepackage{epsfig}
\usepackage[hang,raggedright]{subfigure}

\usepackage{colortbl}
\usepackage{dcolumn}
\newcolumntype{.}{D{.}{.}{1.3}}

\usepackage{txfonts}
\usepackage[hyphens]{url}
\usepackage[subnum]{cases}

\usepackage{eqnarray}
\usepackage{equationarray}
\usepackage{rotating}

\usepackage{tikz}
\usetikzlibrary{positioning}

\newlength{\XWidth}

\newcommand\dashrule{\leavevmode\xleaders\hbox{-}\hfill\kern0pt}

\def\diag{\operatorname{diag}}

\newcommand{\ba}{{\bvec{a}}}
\newcommand{\bb}{\bvec{b}}

\newcommand{\bi}{\bvec{i}}

\newcommand{\bk}{\bvec{k}}
\newcommand{\bl}{\bvec{l}}

\newcommand{\bp}{\bvec{p}}

\newcommand{\bt}{\bvec{t}}
\newcommand{\bu}{\bvec{u}}
\newcommand{\bv}{\bvec{v}}

\newcommand{\bA}{{\bf A}}
\newcommand{\bB}{{\bf B}}

\newcommand{\bF}{{\bf F}}

\newcommand{\bI}{{\bf I}}

\newcommand{\bK}{{\bf K}}

\newcommand{\bM}{{\bf M}}

\newcommand{\bU}{{\bf U}}
\newcommand{\bV}{{\bf V}}
\newcommand{\bW}{{\bf W}}

\newcommand{\bY}{{\bf Y}}

\newcommand{\bsI}{{\boldsymbol I}}

\newcommand{\bsL}{{\boldsymbol L}}

\newcommand{\bsR}{{\boldsymbol R}}

\newcommand{\calL}{{\mathcal{L}}}

\newcommand{\calO}{{\mathcal{O}}}

\newcommand{\blambda}{\mbox{\boldmath $\lambda$}}

\newcommand{\bsigma}{\mbox{\boldmath $\sigma$}}

\newcommand{\bomega}{\mbox{\boldmath $\omega$}}

\newcommand{\bGamma}{\mbox{\boldmath $\Gamma$}}

\newcommand{\bSigma}{\mbox{\boldmath $\Sigma$}}

\newcommand{\bOmega}{\mbox{\boldmath $\Omega$}}

\newcommand{\Real}{\mathbb R}

\newcommand{\1}{\mbox{\boldmath $1$}}

\newcommand{\be}{\begin{eqnarray}}
\newcommand{\ee}{\end{eqnarray}}

\newcommand{\matrixb}{\left[ \begin{array}}
\newcommand{\matrixe}{\end{array} \right]}

\def\*{\circledast}

\newtheorem{theorem}{{Theorem}}[section]
\newtheorem{lemma}{{Lemma}}[section]
\newtheorem{definition}{{Definition}}[section]
\newtheorem{corollary}{{Corollary}}[section]
\newtheorem{remark}{{Remark}}[section]

\newcommand{\bvec}[1]{\boldsymbol{#1}}

\def\vectorize{\operatorname{vec}}
\newcommand{\vtr}[1]{\vectorize\hspace{-.3ex}\left(#1\right)}

\newcommand{\tensor}[1]{\boldsymbol{\mathscr{\MakeUppercase{#1}}}} 
\newcommand{\tA}{\tensor{A}}
\newcommand{\tB}{\tensor{B}}

\newcommand{\tE}{\tensor{E}}
\newcommand{\tF}{\tensor{F}}
\newcommand{\tG}{\tensor{G}}

\newcommand{\tS}{\tensor{S}}

\newcommand{\tX}{\tensor{X}}
\newcommand{\tY}{\tensor{Y}}

{\bfseries}{\itshape}
{\bfseries}{\itshape}

\newcommand{\minitab}[2][l]{\begin{tabular}{@{}#1}#2\end{tabular}}

\usepackage[vlined,ruled,commentsnumbered]{algorithm2e}

\usepackage{aurical}

\DeclareMathAlphabet{\mathscrlow}{OT1}{pzc}{b}{it}

\def\blkdiag{\mathop{\mbox{{\tt{blkdiag}}}}}
\def\bigcircledast{\mathop{\mbox{\fontsize{18}{19}\selectfont $\circledast$}}}

\renewcommand{\bigodot}{\mathop{\mbox{\fontsize{18}{19}\selectfont$\odot$}}}

\usepackage{arydshln}

\usepackage{soul,color}
\setstcolor{red}

\SetKwComment{mtcc}{\% }{}
\SetKwSwitch{Switch}{Case}{Other}{switch}{}{case}{otherwise}{end switch}%
\newcommand{\eqdef}{\stackrel{\triangle}{=}}
\newcounter{example} 

\newenvironment{example}
{\refstepcounter{example}\vspace{10pt}\par\noindent 
\textbf{Example \theexample\ }
\begin{itshape}}
{\end{itshape}}%

\title{CANDECOMP/PARAFAC Decomposition of High-order Tensors Through Tensor Reshaping}

\author{Anh~Huy~Phan$^{*}$, Petr~Tichavsk{\'y} and Andrzej~Cichocki
\thanks{A. H. Phan and A. Cichocki are with the Lab for Advanced Brain Signal Processing, Brain Science Institute, RIKEN, Wakoshi, Japan, e-mail: (phan,cia)@brain.riken.jp.}
\thanks{A. Cichocki is also with System Research Institute, Warsaw, Poland.}
\thanks{P.  Tichavsk{\'y} is with Institute of Information Theory and Automation, Prague, Czech Republic, email: tichavsk@utia.cas.cz.}
\thanks{The work of P. Tichavsk{\'y} was supported by Grant Agency of the Czech Republic 102/09/1278}
}

\begin{document}

\maketitle

\begin{abstract}
%
In general, algorithms for order-3 CANDECOMP/\-PARAFAC (CP), also coined canonical polyadic decomposition (CPD), are easily to implement and can be extended to higher order CPD.
Unfortunately, the algorithms become computationally demanding, and they are often not applicable to higher order and relatively large scale tensors.
In this paper, by exploiting the uniqueness of CPD and the relation of a tensor in Kruskal form and its unfolded tensor, we propose a fast approach to deal with this problem.
Instead of directly factorizing the high order data tensor, the method decomposes an unfolded tensor with lower order, e.g., order-3 tensor. 
On basis of the order-3 estimated tensor, a structured Kruskal tensor of the same dimension as the data tensor is then generated, and decomposed to find the final solution using fast algorithms for the structured CPD.
In addition, strategies to unfold tensors are suggested and practically verified in the paper.
\end{abstract}



\section{Introduction}
\label{sec:introduction}

CANDECOMP/\-PARAFAC \cite{Harshman,Carroll_Chang}, also known as Canonical polyadic decomposition (CPD), 
is a common tensor factorization which has found applications such as in chemometrics\cite{Andersen03,bro_book,Tomasithesis}, telecommunication\cite{DC07,Sidiropoulos2000}, analysis of fMRI data\cite{Andersen}, time-varying EEG spectrum \cite{Morup,DBLP:journals/sigpro/BeckerCAHM12}, data mining \cite{springerlink:10.1007,ACM-TKDD-Link-Prediction},  separated representations for generic functions involved in quantum mechanics or kinetic theory descriptions of materials \cite{NME:NME2710}, classification, clustering \cite{Shashua}, compression \cite{NLA:NLA297,NLA:NLA344,NLA:NLA765}. 
Although the original decomposition and applications were developed for three-way data, the model was later widely extended to higher order tensors.
For example, P. G. Constantine \emph{et al.} \cite{Constantine} modeled the pressure measurements along the combustion chamber as order-6 tensors corresponding to the flight conditions - Mach number, altitude and angle of attack,  and the wall temperatures in the combustor and the turbulence mode. 
Hackbusch, Khoromskij, and Tyrtyshnikov \cite{Hackbusch03hierarchicalkronecker} and Hackbusch and Khoromskij \cite{Hackbusch2007697} investigated CP approximation to operators and functions in high dimensions.  Oseledets and Tyrtyshnikov \cite{Oseledets:2009p9706} approximated the Laplace operator and the general second-order operator which appears in the Black-Scholes equation for multi-asset modeling to tackle the dimensions up to $N$ = 200. 
In neuroscience, M. M{\o}rup \emph{et al.} \cite{Morup} analyzed order-4 data constructed from EEG signals in the time-frequency domain.
Order-5 tensors consisting of  dictionaries $\times$ timeframes $\times$ frequency bins $\times$ channels $\times$ trials-subjects \cite{Phan-CIA_NOLTA10} built up from EEG signals were shown to give high performance in BCI based on EEG motor imagery.
In object recognition (digits, faces, natural images), CPD was used to extract features from order-5 Gabor tensors including hight $\times$  width  $\times$ orientation $\times$ scale  $\times$ images \cite{Phan-CIA_NOLTA10}.
 
%
%


In general, many CP algorithms for order-3 tensor can be straightforwardly extended to decompose higher order tensors.
For example, there are numerous algorithms for CPD including the alternating least squares (ALS) algorithm \cite{Carroll_Chang,Harshman} with line search extrapolation methods \cite{Harshman,Nwaytoolbox,Rajih-SIAM,Tomasithesis,DBLP:journals/tsp/ChenHQ11}, rotation \cite{Paaterorotation10} and compression \cite{Kiers98}, or all-at-once algorithms such as the OPT algorithm \cite{JChem-CPOPT}, the conjugate gradient algorithm for nonnegative CP, the PMF3, damped Gauss-Newton (dGN) algorithms \cite{Paatero97,Tomasithesis} and fast dGN \cite{Petr_10,DBLP:journals/corr/abs-1205-2584,Phan_LMplus}, or algorithms based on joint diagonalization problem \cite{deLathauwer-JAD,RH:08,DBLP:journals/tsp/LathauwerC08}.  
The fact is that the algorithms become more complicated, computationally demanding, and often not applicable to relatively large scale tensors. For example, complexity of gradients of the cost function with respect to factors grows linearly with the number of dimensions $N$.  It has a computational cost of order $\mathcal O\left(\displaystyle N R \prod_{n=1}^{N} I_n\right)$ for a tensor of size $I_1 \times I_2 \times \cdots \times I_N$. More tensor unfoldings $\bY_{(n)}$ ($n = 2, 3, \ldots, N-1$) means more time consuming due to accessing non-contiguous blocks of data entries and shuffling their orders in a computer. 
In addition, line search extrapolation methods\cite{Harshman,bro_book,Nwaytoolbox,Francthesis,Tomasithesis,Rajih-SIAM} become more complicated, and demand high computational cost to build up and solve $(2N-1)$- order polynomials.
The rotation method \cite{Paaterorotation10} needs to estimate $N$ rotation matrices of size $R \times R$ with a whole complexity per iteration of order $\calO(N^3R^6)$.

%
%
%



Recently, a Cram\'er-Rao Induced Bounds (CRIB) on attainable squared angular error of factors in the CP decomposition has been proposed in \cite{2012arXiv1209.3215T}. The bound is valid under the assumption that the decomposed
tensor is corrupted by additive Gaussian noise which is
independently added to each tensor element. In this paper we use
the results of \cite{2012arXiv1209.3215T} to design the tensor unfolding strategy which ensures as little deterioration of accuracy as possible.
This strategy is then verified in the
simulations.

By exploiting the uniqueness of CPD under mild conditions and the relation of a tensor in the Kruskal form \cite{Kolda08} and its unfolded tensor, we propose a fast approach for high order and relatively large-scale  CPD. Instead of directly factorizing the high order data tensor, the approach decomposes an unfolded tensor in lower order, e.g., order-3 tensor.  A structured Kruskal tensor of the same dimension of the data tensor is then generated, and decomposed to find the desired factor matrices. We also proposed the fast ALS algorithm to factorize the structured Kruskal tensor. 

The paper is organized as follows. Notation and the CANDECOMP/PARAFAC are briefly reviewed in Section~\ref{sec:notation}.
The simplified version of the proposed algorithm is presented in Section~\ref{sec::lowrank3w2Nw}.
Loss of accuracy is investigated in Section~\ref{sec::CRIBloss}, and an efficient strategy for tensor unfolding is summarized in Section~\ref{sec::unfoldingstrategy}.
For difficult scenario decomposition, we proposed a new algorithm in Section~\ref{sec::3wtoNw-lowrank}. 
Simulations are performed on random tensors and real-world dataset in Section~\ref{sec::simulation}.
Section \ref{sec:conclusion} concludes the paper.

%
%
%
%
%
%
%
%
%
%
%
%
\section{CANDECOMP/PARAFAC (CP) decomposition}\label{sec:notation}

Throughout the paper, we shall denote tensors by bold calligraphic letters, e.g., $\tA \in \Real^{I_1 \times I_2 \times \cdots \times I_N}$,
matrices by bold capital letters, e.g., $\bA$ =$[\ba_1,\ba_2, \ldots, \ba_R] \in \Real^{I \times R}$, and vectors by bold italic letters, e.g., $\ba_j$ or $\bsI = [I_1, I_2,\ldots, I_N]$.
%
A vector of integer numbers is denoted by colon notation such as $\bk = i\text{:}j = [i, i+1, \ldots, j-1, j]$. For example, we denote $1\text{:}n= [1, 2, \ldots, n]$.
The Kronecker product, the Khatri-Rao (column-wise Kronecker) product, and the (element-wise) Hadamard product  are  denoted respectively by  $\otimes, \odot, \circledast$ \cite{Kolda08,NMF-book}.

\begin{definition}\label{def_Kruskal}{\bf(Kruskal form (tensor) \cite{Kolda08,Loanslecture5})} 
A tensor $\tX \in \Real^{I_1 \times I_2 \times \cdots \times I_N}$ is in Kruskal form if
\be
{\tX}  &=&  \sum\limits_{r = 1}^R  \lambda_r \, {\ba^{(1)}_{r}  \circ \ba^{(2)}_{r} \circ  \cdots  \circ \ba^{(N)}_{r}}  ,  \\
&\eqdef& \llbracket \blambda;  \bA^{(1)}, \bA^{(2)}, \ldots, \bA^{(N)} \rrbracket, \quad \blambda  = [\lambda_1, \lambda_2, \ldots, \lambda_R].
\label{equ_CP_nolambda}
\ee
where symbol ``$\circ$"  denotes the outer product,  $\bA^{(n)}=[\ba^{(n)}_1, \ba^{(n)}_2,\ldots,\ba^{(n)}_R]$ $ \in \Real ^{I_n \times R}$, $(n=1,2, \ldots, N)$ are factor matrices, 
$\ba^{(n) T}_r \ba^{(n)}_r = 1$, for all $r$ and $n$, and $\lambda_1 \ge \lambda_2 \ge \cdots \ge \lambda_R >0$.\end{definition}

\begin{definition}\label{def_CP}{\bf(CANDECOMP/PARAFAC (CP)\cite{Hitchcock1927,Harshman,Carroll_Chang,Loanslecture5} )}
Approximation of order-$N$ data tensor ${\tY} \in \Real^{I_1 \times I_2 \times \cdots \times I_N}$ by a rank-$R$ tensor in the Kruskal form means
\be
{\tY}  &=&  \tensor{\widehat Y} +\tE,  
\ee
where $\tensor{\widehat Y} = \llbracket \blambda;  \bA^{(1)}, \bA^{(2)}, \ldots, \bA^{(N)} \rrbracket$,  so that  $\|\tY - \tensor{\widehat Y}\|_F^2$ is minimized.
 
\end{definition}

There are numerous algorithms for CPD including alternating least squares (ALS) or all-at-once optimization algorithms, or based on joint diagonalization. In general, most CP algorithms which factorize order-$N$ tensor often face high computational cost due to computing gradients and (approximate) Hessian, line search and rotation. Table~\ref{tab_complexity} summarizes complexities of major computations in popular CPD algorithms.
Complexity per iteration of a CP algorithm can be roughly computed based on Table~\ref{tab_complexity}. For example, the ALS algorithm with line search has a complexity of order $\calO(NRJ + 2^N RJ + NR^3) = \calO(2^N RJ + NR^3)$.

\begin{table}
\centering
\caption{Complexities per iteration of major computations in CPD algorithms. 
$J = \prod_{n = 1}^{N} I_n$, $T =  \sum_{n=1}^{N} I_n$.}\label{tab_complexity}
\begin{tabular}{ll}
\multicolumn{1}{c}{Computing Process}  & \multicolumn{1}{c}{Complexity} \\\hline
Gradient\cite{Carroll_Chang,Harshman}& $\calO\left(NR J\right)$ \\
Fast gradient\cite{DBLP:journals/corr/abs-1204-1586}   & $\calO\left(R J \right)$ \\ 
(Approximate) Hessian and its inverse\cite{Paatero97,Tomasithesis} & $\calO\left(R^3 T^3\right)$ \\
Fast (approximate) Hessian and its inverse\cite{DBLP:journals/corr/abs-1205-2584,2012arXiv1209.3215T}  & $\calO\left(R^2T + N^3R^6\right)$ \\
Exact line search\cite{Harshman,bro_book,Tomasithesis} & $\calO\left(2^N R J  \right)$\\
Rotation\cite{Paaterorotation10}  &   $\calO\left(N^3 R^6 \right)$\\
\end{tabular}
\end{table}





\section{CPD of unfolded tensors}\label{sec::lowrank3w2Nw}

In order to deal with existing problems for high order and relatively large scale CPD, the following process is proposed:
\begin{enumerate}
\item Reduce the number of dimensions of the tensor $\tY$ to a lower order (e.g., order-3) through tensor unfolding $\tY_{\llbracket \bl \rrbracket}$ which is defined later in this section. 
\item Approximate the unfolded tensor $\tY_{\llbracket \bl \rrbracket}$ by an order-3 tensor ${\tensor{\widehat Y}}_{\llbracket \bl \rrbracket}$ in the Kruskal form. 
Dimensions of $\tY_{\llbracket \bl \rrbracket}$ which are relatively larger than rank $R$ can be reduced to $R$ by the Tucker compression\cite{Tucker66,tucker64extension,Lathauwer_HOOI,Andersson98} 
prior to CPD although it is not a lossless compression.
In such case, we only need to decompose an $R \times R \times R$ dimensional tensor.
\item Estimate the desired components of the original tensor $\tY$ on basis of the tensor ${\tensor{\widehat Y}}_{\llbracket \bl \rrbracket}$ in the Kruskal form. 
\end{enumerate}

The method is based on an observation that unfolding of a Kruskal tensor also yields a Kruskal tensor. Moreover due to uniqueness of CPD under ``mild'' conditions, the estimated components along the unfolded modes are often good approximates to components for the full tensor.
In the sequel, we introduce basic concepts that will be used in the rest of this paper.
Loss of accuracy in decomposition of the unfolded tensors is analyzed theoretically based on the CRIB.

\begin{definition}[\bf{Reshaping}]\label{def_reshape}
The reshape operator for a tensor $\tY \in \Real^{I_1 \times I_2 \times \cdots \times I_N}$
to a size specified by a vector $\bsL = [L_1, L_2, \ldots, L_M]$ with $\prod_{m = 1}^{M} L_m = \prod_{n = 1}^{N} I_n$ returns an order-$M$ tensor $\tX$, such that
$ \vtr{\tY} = \vtr{\tX}$, and is expressed as
$
    \tX = {\tt{reshape}}(\tY, \bsL)  
    \;  \in \Real^{L_1 \times L_2 \times \cdots \times L_M} $.
\end{definition}


\begin{definition}
[{\bf Tensor transposition} \cite{DBLP:journals/siammax/RagnarssonL12}]
If $\tA \in \Real^{I_1 \times \cdots \times I_N}$ and $\bp$ is a permutation of $[1, 2, \ldots, N]$, then
$\tA^{<\bp>} \in \Real^{I_{p_1} \times \cdots \times  I_{p_N}} $
denotes the $\bp$-transpose of $\tA$ and is defined by
\be
	\tA^{<\bp>}(i_{p_1} , \ldots, i_{p_N}) = \tA(i_1, \ldots, i_N),  \quad  \1 \leq \bi \leq \bsI = [I_1, I_2, \ldots, I_N]. \label{equ_tensor_transpose}
\ee
\end{definition}

\begin{definition}
[{\bf Generalized tensor unfolding}]\label{def_unfolding}
Reshaping a $\bp$-transpose $\tY^{<\bp>}$ to an order-$M$ tensor of size $\bsL = [L_1, L_2, \ldots, L_M]$ with $\displaystyle L_m = \prod_{k \in \bl_m} I_k$, where $[\bl_1, \bl_2, \ldots, \bl_M ] \equiv [p_1, p_2, \ldots, p_N]$, $\bl_m = [\bl_m(1), \ldots, \bl_m(K_m)]$ 
\be
	\tY_{\llbracket \bl \rrbracket}   \eqdef
	{\tt{reshape}}(\tY^{<\bp>}, \bsL), \qquad \bl = [\bl_1, \bl_2, \ldots, \bl_M ] .
\ee
\end{definition}

\begin{remark}\hspace{0ex}\\[-2em]
\begin{enumerate}
\item If $\bl = [n, (1\textup{:}n-1, n+1\textup{:}N)]$, then $\tY_{\llbracket \bl \rrbracket} = \bY_{(n)}$ is mode-$n$ unfolding.
\item If  $\tY$ is an order-4 tensor, then $\tY_{\llbracket 1, 2, (3,4) \rrbracket}$ is an order-3 tensor of size $I_1 \times I_2 \times I_3I_4$.
\item If  $\tY$ is an order-6 tensor, then $\tY_{\llbracket (1, 4), (2,5), (3,6) \rrbracket}$ is an order-3 tensor of dimension $I_1 I_4 \times I_2I_5 \times I_3I_6$.
\end{enumerate}
\end{remark}

We denote Khatri-Rao product of a set of matrices $\bU^{(n)}$, $n = 1, 2, \ldots, N$, as $\displaystyle \bigodot_{n = 1}^{N} \bU^{(n)} \eqdef \bU^{(N)} \odot \bU^{(N-1)} \odot \cdots \odot \bU^{(1)}$.

\begin{lemma}\label{lem_CPfolding}
Unfolding of a rank-$R$ tensor in the Kruskal form $\tY = \llbracket \blambda; \bA^{(1)} , \bA^{(2)} , \ldots,  \bA^{(N)} \rrbracket$ returns an order-$M$ rank-$R$ Kruskal tensor $\tY_{\llbracket \bl \rrbracket}$, $\bl = [\bl_1, \bl_2, \ldots, \bl_M]$, given by
\be
	\tY_{\llbracket \bl \rrbracket} = \llbracket \blambda; \bB^{(1)} , \bB^{(2)} , \ldots,  \bB^{(M)} \rrbracket ,
\ee
where
$
	\bB^{(m)} = \bigodot_{k \in \bl_m} \bA^{(k)} 
	\;  \in \Real^{\scriptstyle \left(\prod_{k \in \bl_m} I_k\right) \times R}$ $(m = 1, 2, \ldots, M)$ are merging factor matrices.
\end{lemma}

\begin{remark}\hspace{0ex}\\[-2em]
\begin{enumerate}
\item If $\bl = [n, (1\textup{:}n-1, n+1\textup{:}N)]$, then 
$
	\tY_{\llbracket \bl \rrbracket} = \tY_{(n)} = \bA^{(n)} \, \diag(\blambda) \left(\bigodot_{k \neq n} \bA^{(k)} \right)^T
$.
\item If $\bl = [(1\textup{:}n), (n+1\textup{:}N)]$, then 
$
	\tY_{\llbracket \bl \rrbracket} =  \left(\bigodot_{k = 1}^{n} \bA^{(k)} \right) \, \diag(\blambda) \left(\bigodot_{k = n+1}^{N} \bA^{(k)} \right)^T
$.
\item For an order-4 Kruskal tensor $\tY$, $\tY_{\llbracket 1, 2, (3,4) \rrbracket} = \llbracket \blambda; \bA^{(1)} , \bA^{(2)} , \bA^{(4)}  \odot \bA^{(3)} \rrbracket$.
\end{enumerate}
\end{remark}

\begin{corollary}\label{corollary_CPrnk1}
An order-$K$ tensor $\tB_{m_r}$ of size $I_{l_{m1}} \times I_{l_{m2}} \times \cdots \times I_{l_{mK}}$, $\bl_m = [{l_{m1}} \ldots {l_{mK}}]$ folded from the $r$-th column vector $\bb^{(m)}_r$ of $\bB^{(m)}$, i.e., $\vtr{\tB_{m_r}} = \bb^{(m)}_r$ is a rank-1 tensor 
\be
	\tB_{m_r} =   \ba^{(l_{m1})}_r \circ \ba^{(l_{m2})}_r \circ \cdots \circ \ba^{(l_{mK})}.
\ee
\end{corollary}

In practice for real data, folded tensors $\tB_{m_r}$ are not exact rank-1 tensors but can be approximated by rank-1 tensors composed from components corresponding modes in $\bl_m$. In other words, computing the leading-left singular vector of the mode-$k$ unfolding $\left[\tB_{m_r}\right]_{(k)}$ is the simplest approach to recover $\ba^{(l_{mk})}_r$ from $\bb^{(m)}_r$ for $k = 1, 2, \ldots, K$. Pseudo-code of this simple algorithm for unFolding CPD (FCP) is described in Algorithm~\ref{alg_3w2Nw_rank1}.
The more complex and efficient algorithm is discussed later.

\def\baselinestretch{1} 
\setlength{\algomargin}{1em}
\begin{algorithm}[t!]
\SetFillComment
\SetSideCommentRight
\caption{{\tt{rank-one FCP}}\label{alg_3w2Nw_rank1}}
\DontPrintSemicolon \SetFillComment \SetSideCommentRight
\KwIn{Data tensor $\tY$:  $(I_1 \times I_2 \times \cdots \times I_N)$, rank $R$, \linebreak 
Unfolding rule $\bl = [\bl_1, \bl_2, \ldots, \bl_M]$ where $\bl_m = [ l_m(1), \ldots, l_m(K_m)]$} 
\KwOut{$\blambda \in \Real^{N}$, $N$ matrices $\bA^{(n)} \in \Real^{I_n \times R}$} \SetKwFunction{mreshape}{reshape}
\SetKwFunction{permute}{permute} 
\SetKwFunction{svds}{svds} 
\SetKwFunction{eigs}{eigs} 
\SetKwFunction{rankone}{rank1CP} 
\SetKwFunction{lowrank}{lowrankCP} 
\SetKwFunction{CPD}{CPD} 
\SetKwFunction{psCPD}{structuredCPD} 
\SetKwFunction{MwtoNw}{FCP} 
\SetKwFunction{tucker}{TD} 
\Begin{
 {\mtcc{Stage 1: Tensor unfolding and optional compression \dashrule}}
\nl $ \llbracket \tG, \bU^{(1)},  \ldots, \bU^{(M)} \rrbracket = \tucker(\tY_{\llbracket \bl \rrbracket},\min(\bI, R))${\mtcc*{Tucker decomposition of order-$M$ $\tY_{\llbracket \bl \rrbracket}$}}
\;\vspace{-1ex}
{\mtcc{Stage 2: CPD of the unfolded (and reduced) tensor \dashrule}}
\nl $ \llbracket \blambda; \bB^{(1)},  \ldots, \bB^{(M)} \rrbracket = \CPD(\tG,R)${\mtcc*{order-$M$ CPD of the core tensor }}
\nl \lFor{$m  = 1, 2, \ldots, M$}{
$\bB^{(m)} \leftarrow \bU^{(m)} \, \bB^{(m)}${\mtcc*{Back projection of TD}}
}
\;\vspace{-1ex}
 {\mtcc{Stage 3: Rank-one approximation to merging components \dashrule}}
\nl \For{$m  = 1, 2, \ldots, M$}{
	\nl \For{$r  = 1, 2, \ldots, R$}{
	\nl $ \llbracket g; \ba^{(l_{m1})}_r,  \ldots, \ba^{(l_{mK})}_r \rrbracket = 
		\tucker({\tt{reshape}}(\bb^{(m)}_r,{[I_{l_{m}(1)}, \ldots, I_{l_{m}(K_m)}]}),1)$\;
	\nl $\lambda_r \leftarrow \lambda_r  \, g$	
	}
}
}
\BlankLine
\minitab[p{.95\linewidth}]{\hline
$\tucker(\tY, \bsR)$:  rank-$\bsR$ Tucker decomposition of order-$N$ tensor $\tY$ where $\bsR = [R_1, R_2, \ldots, R_N]$.\\
${\tensor{\widehat{Y}}}$ = \CPD($\tY$, $R$, $\tY_{init}$):  approximates an order-$N$ tensor or a tensor in the Kruskal form $\tY$ by a rank-$R$ Kruskal tensor ${\tensor{\widehat{Y}}}$ using initial values $\tY_{init}$.
}
\end{algorithm}
\def\baselinestretch{1.5}


\subsection{Selecting an unfolding strategy}\label{sec::CRIBloss}

For (noiseless) tensors which have exact rank-$R$ CP decompositions without (nearly) collinear components, factors computed from unfolded tensors can be close to the true solutions. However, for real data tensor, there exists loss of accuracy when using the rank-one approximation approach. The loss can be affected by the unfolding, or by the rank-$R$ of the decomposition, especially when $R$ is under the true rank of the data tensor.
This section analyzes such loss based on comparing CRIBs on the first component $\ba^{(1)}_1$ of CPDs of the full tensor and its unfolded version. We use $\ba_1$ a shorthand notation for $\ba^{(1)}_1$.
The results of this section give us an insight into how to unfold a tensor without or possibly minimal loss of accuracy.

The accuracy loss in decomposition of unfolded tensor is defined as the loss of 
CRIB \cite{DBLP:conf/icassp/TichavskyK11,ZbynekSSP11,2012arXiv1209.3215T}
on components of the unfolded tensor through the unfolding rule $\bl$ compared with CRIB on components of the original tensor.
%
For simplicity, we consider tensors in the Kruskal form (\ref{def_Kruskal}) which have $ \ba^{(n)T}_r \, \ba^{(n)}_s = c_n$, for all $n$, $r \neq s$, and $\ba^{(n)T}_r \, \ba^{(n)}_r = 1$, $-1 \le c_n \le1$. Coefficients $c_n$ are called degree of collinearity.

\subsubsection{Loss in unfolding order-4 tensors}\label{sec::loss4D}

For order-4 tensors,  since $\mbox{CRIB}(\ba_1)$ is largely independent of $c_1$ unless $c_1$ is close to $\pm 1$ \cite{2012arXiv1209.3215T}, we consider the case when $c_1=0$. Put $h=c_2c_3c_4$ and $\theta = \frac{\sigma^2}{\lambda_1^2}$. From \cite{2012arXiv1209.3215T} and Appendix~\ref{sec::CRIB}, we get CRIBs for rank-2 CPD in explicit forms:
\begin{eqnarray}
\mbox{CRIB}(\ba_1)&=&\frac{\theta}{1-h^2}\left[I_1-1+\frac{c_2^2c_3^2+c_2^2c_4^2+c_3^2c_4^2-3h^2}
{1+2h^2-c_2^2c_3^2-c_2^2c_4^2-c_3^2c_4^2}\right], \\
\mbox{CRIB}_{\llbracket 1,2,(3,4)\rrbracket}(\ba_1)&=&\frac{\theta}{1-h^2}\left[I_1-3+\frac{1}{1-c_2^2} + \frac{1}{1-c_3^2c_4^2}\right], \label{eq_34}\\
\mbox{CRIB}_{\llbracket 1,(2,3),4)\rrbracket}(\ba_1)&=&\frac{\theta}{1-h^2}\left[I_1 -3+\frac{1}{1-c_2^2c_3^2} + \frac{1}{1-c_4^2}\right]. \label{eq_12}
\end{eqnarray}
In general, $\mbox{CRIB}(\ba_1)\leq\mbox{CRIB}_{\llbracket1,2,(3,4)\rrbracket}(\ba_1)$. The equality is achieved for $c_2=0$.
\begin{eqnarray}
\mbox{CRIB}(\ba_1)_{c_2=0}&=&
\mbox{CRIB}_{\llbracket1,2,(3,4)\rrbracket}(\ba_1)_{c_2=0}= \theta \left(I_1-2+\frac{1}
{1-c_3^2c_4^2}\right). \label{eq_ortho4}
\end{eqnarray}
It means that if modes 1 and 2 comprise (nearly) orthogonal components, the tensor unfolding $[1,2,(3,4)]$ does not affect the accuracy of the decomposition.

From (\ref{eq_34}) and (\ref{eq_12}), it is obvious that  $\mbox{CRIB}_{\llbracket 1,2,(3,4)\rrbracket}(\ba_1)   \le  \mbox{CRIB}_{\llbracket 1,(2,3),4)\rrbracket}(\ba_1)$ if $c_2^2 \le c_4^2$.
This indicates that {\textit{\textbf{collinearities of modes to be unfolded should be  higher than those of other modes in order to reduce the loss of accuracy in estimating $\ba_1$.}}} 
Note that the new factor matrices yielded through tensor unfolding have lower collinearity than those of original ones.
Moreover, tensors with high collinear components  are always more difficult to decompose than ones with lower collinearity \cite{MITCHELL94,CEM:CEM1236,Paatero97}.
Hence, it is natural to unfold modes with highest collinearity so that the CPD becomes easier.
This rule also holds for higher rank $R$, and is illustrated in a particular case when $c_1 = c_3 = 0$
\be
\mbox{CRIB}_{\llbracket1,2,(3,4)\rrbracket}(\ba_1) &=& \theta \left( I_1 - R + \frac{R-1}{1- c_2^2} \right),  \\
\mbox{CRIB}_{\llbracket1,(2,3),4\rrbracket}(\ba_1) &=& \theta \left(I_1 - R + \frac{R-1}{1- c_4^2}\right) .
\ee
The unfolding $[1,2,(3,4)]$ is more efficient than the unfolding $[1,(2,3),4]$ when $|c_2| < |c_4|$, 
although this unfolding still causes some loss of CRIB despite of $c_4$ since
\be
	\mbox{CRIB}(\ba_1)=   \theta\left(I_1-R + \frac{R-1}{1 - c_2^2 c_4^2} \right) \le 
\theta\left( I_1-R+\frac{R-1}
{1-c_2^2}\right) = \mbox{CRIB}_{\llbracket 1,2,(3,4)\rrbracket}(\ba_1), \;  \text{ for all }  c_2. \notag
\ee
Moreover the loss is significant when $c_4$ is small enough.
Note that for this case, the unfolding $[1,3, (2,4)]$ is suggested because it does not cause any loss according to the previous rule.

In other words, {\textit{\textbf{modes which comprise orthogonal or low collinear components (i.e., $c_n\approx 0$) should not fold with the other modes, unless the other modes have nearly orthogonal columns as well.}}}

%

\begin{figure}[t!]
\psfrag{I2}[t][t]{\scalebox{1}{\color[rgb]{0,0,0}\setlength{\tabcolsep}{0pt}\begin{tabular}{c}\\[-2.5em]\small $I$\end{tabular}}}%
\centering
\subfigure[($c_1,c_2,c_3,c_4$) = ({\color{red}{0.1,0.1}},0.1,0.1).]
{
\includegraphics[width=.42\linewidth, trim = 0.0cm .5cm 0cm 1cm,clip=true]{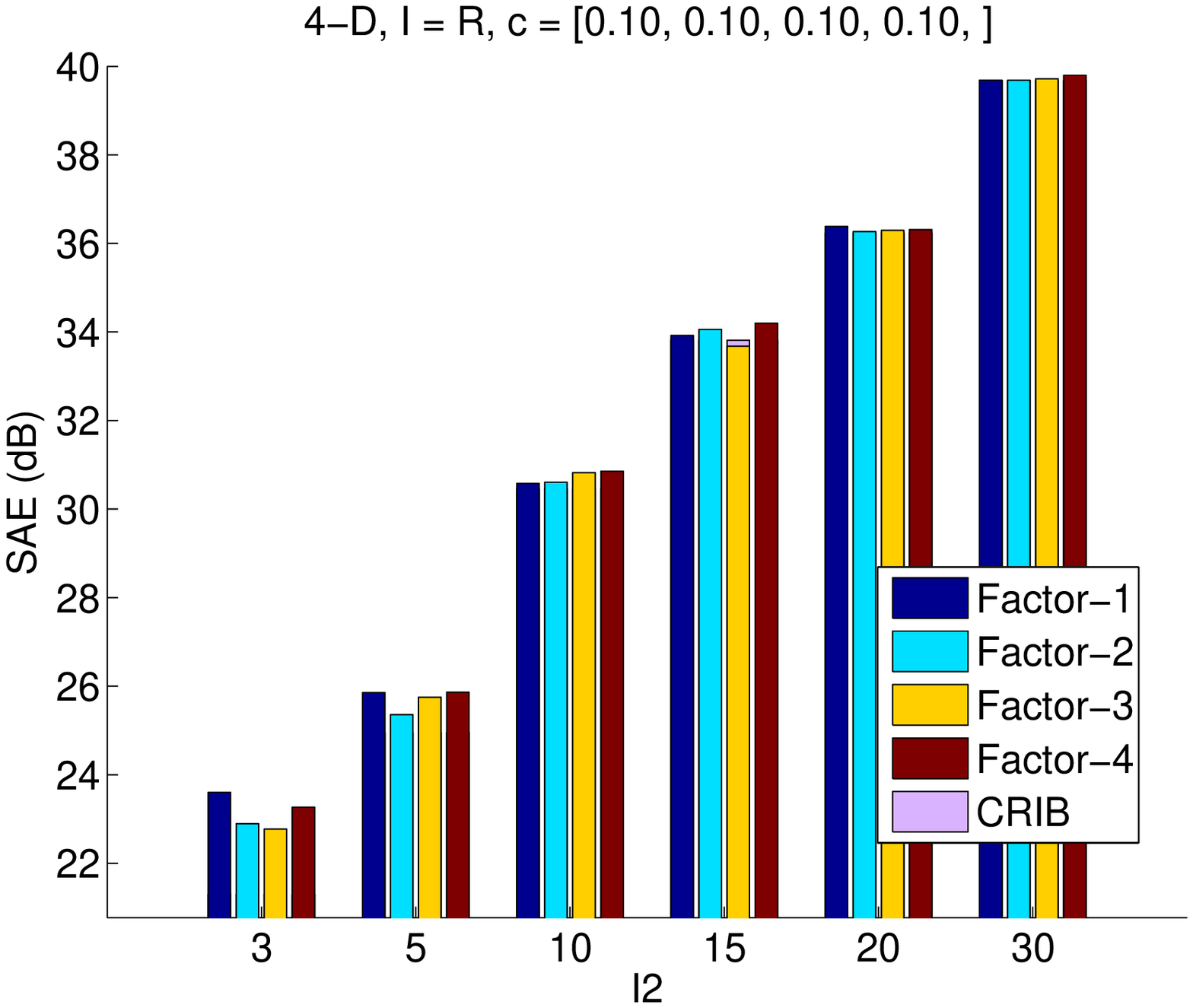}\label{fig_N4_c1111_msaevscrib_allcomps}}
\hfill
\subfigure[($c_1,c_2,c_3,c_4$) = ({\color{red}{0.1,0.1}},0.9,0.1).]
{
\includegraphics[width=.42\linewidth, trim = 0.0cm 0.5cm 0cm 1cm,clip=true]{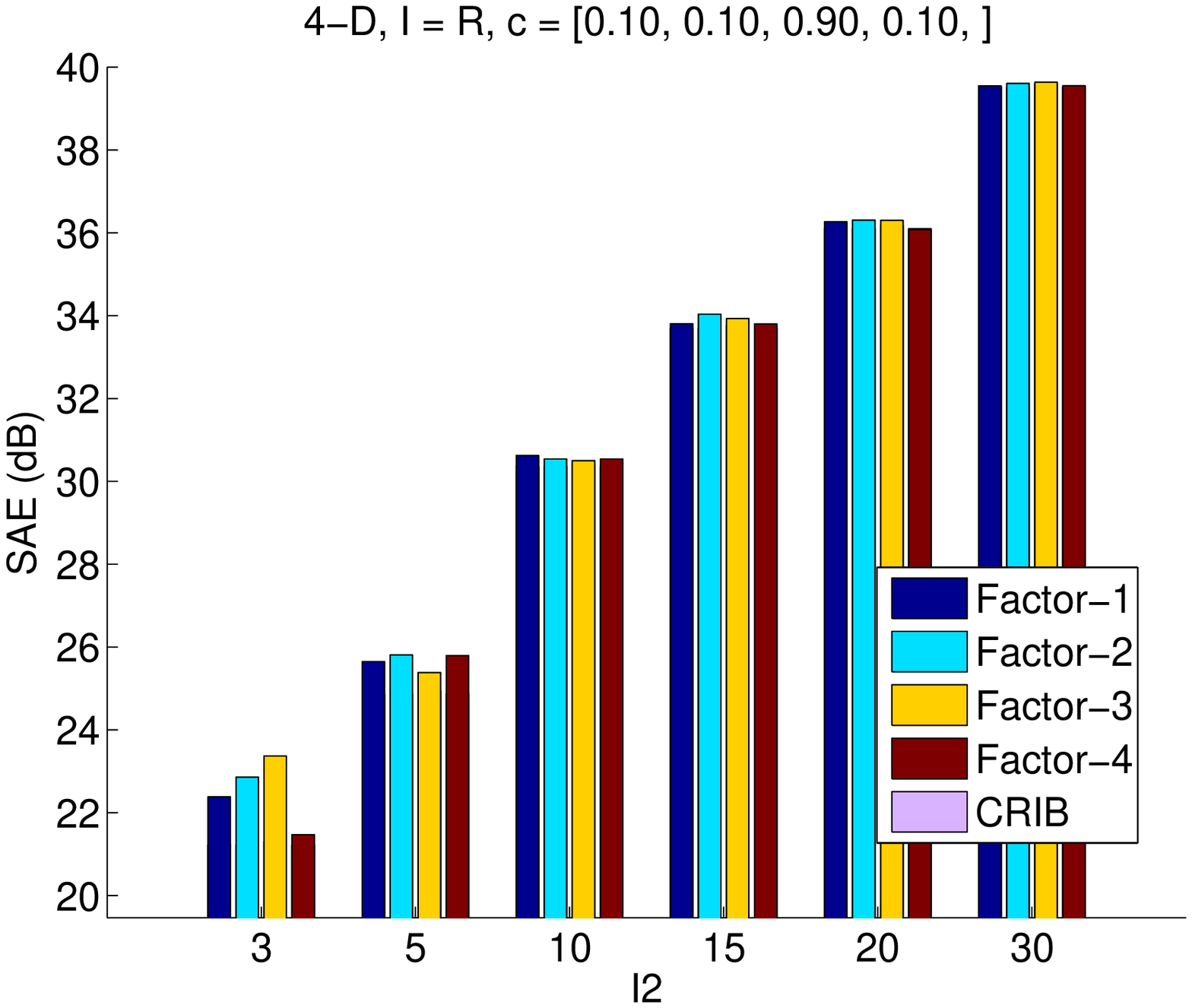}\label{fig_N4_c1191_msaevscrib_allcomps}}
\hfill
\subfigure[($c_1,c_2,c_3,c_4$) = ({\color{red}{0.1,0.1}},0.9,0.9).]
{
\includegraphics[width=.42\linewidth, trim = 0.0cm 0.5cm 0cm 1cm,clip=true]{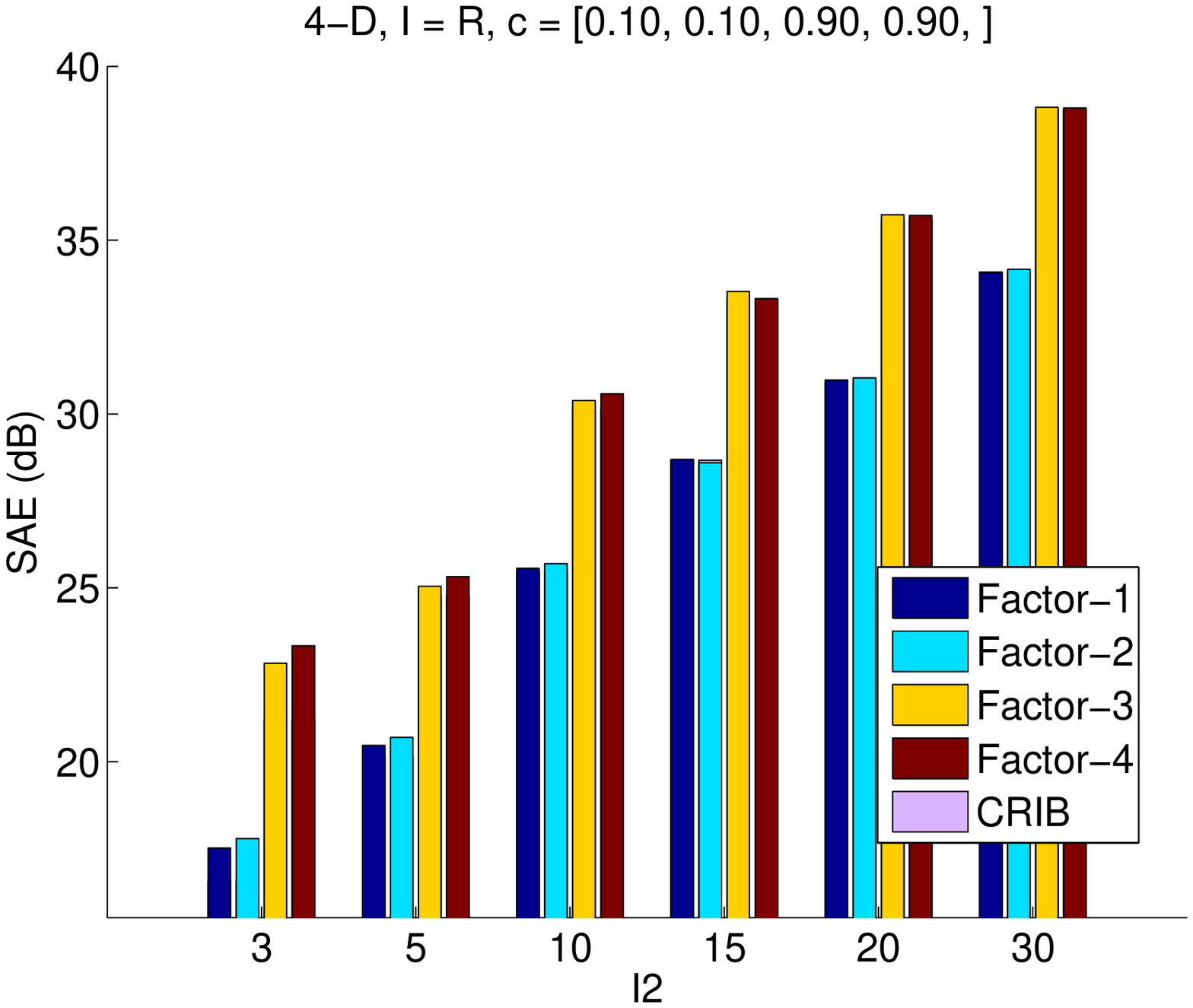}\label{fig_N4_c1199_msaevscrib_allcomps}}
\hfill
\subfigure[($c_1,c_2,c_3,c_4$) = (0.1,0.9,{\color{red}{0.1}},0.9).]
{
\includegraphics[width=.42\linewidth, trim = 0.0cm 0.5cm 0cm 1cm,clip=true]{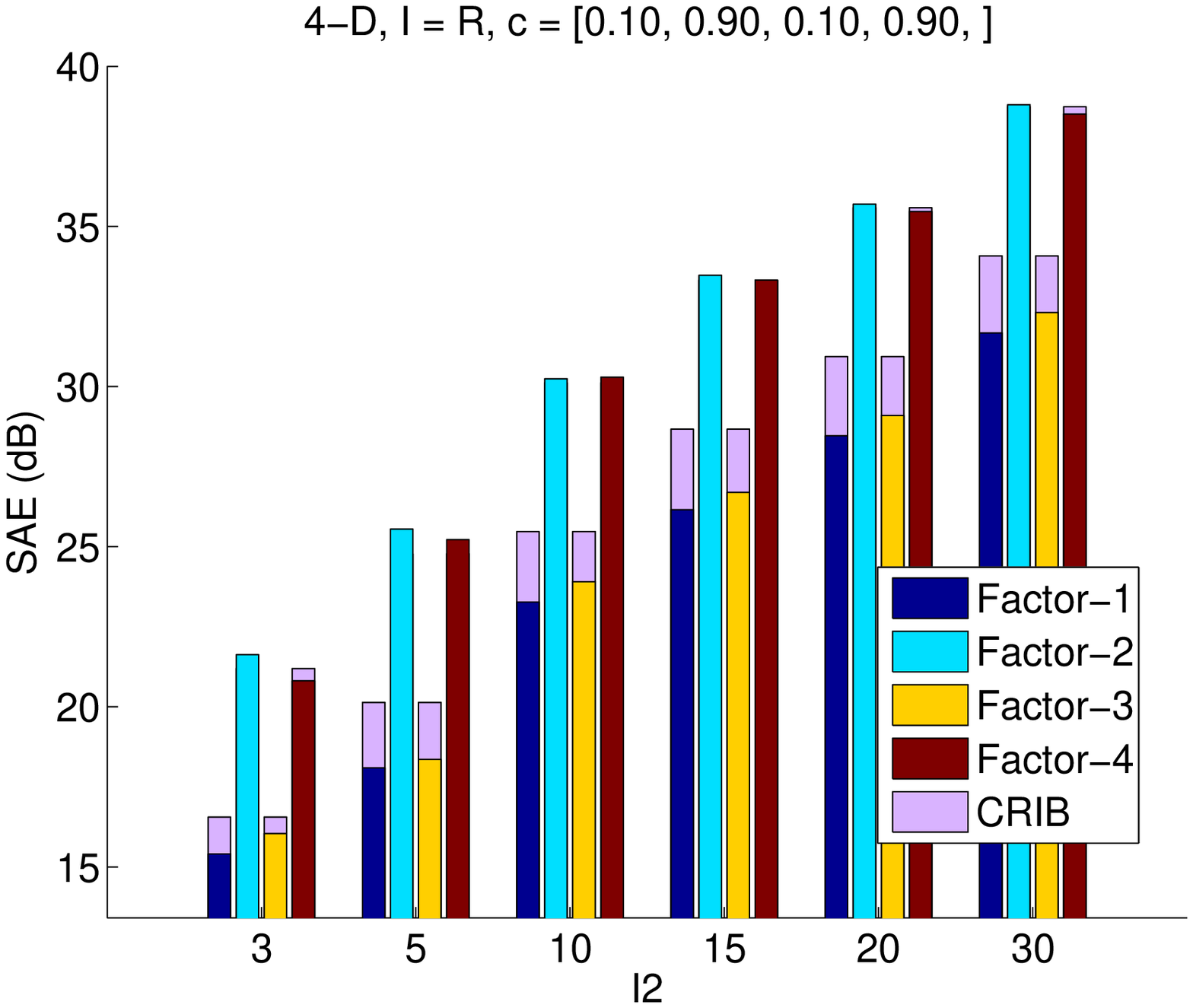}\label{fig_N4_c1919_msaevscrib_allcomps}}
%
\hfill
\subfigure[($c_1,c_2,c_3,c_4$) = (0.1,0.9,{\color{red}{0.9,0.9}}).]
{
\includegraphics[width=.42\linewidth, trim = 0.0cm 0.5cm 0cm 1cm,clip=true]{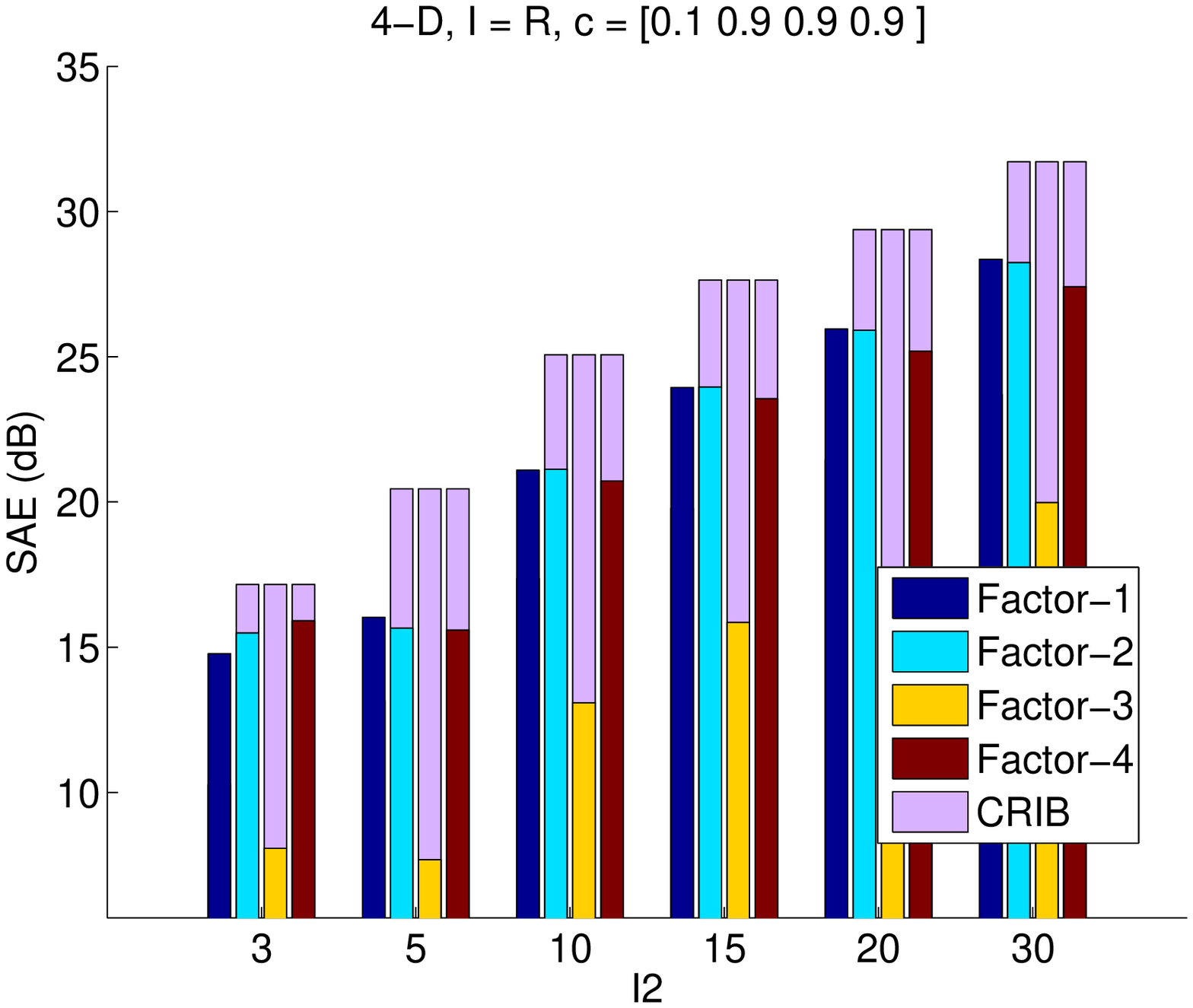}\label{fig_N4_c1999_msaevscrib_allcomps}}
%
%
%
\hfill
\subfigure[($c_1,c_2,c_3,c_4$) = (0.9,0.9,0.9,0.9), SNR = 30 dB.]
{
\includegraphics[width=.42\linewidth, trim = 0.0cm 0.5cm 0cm 1cm,clip=true]{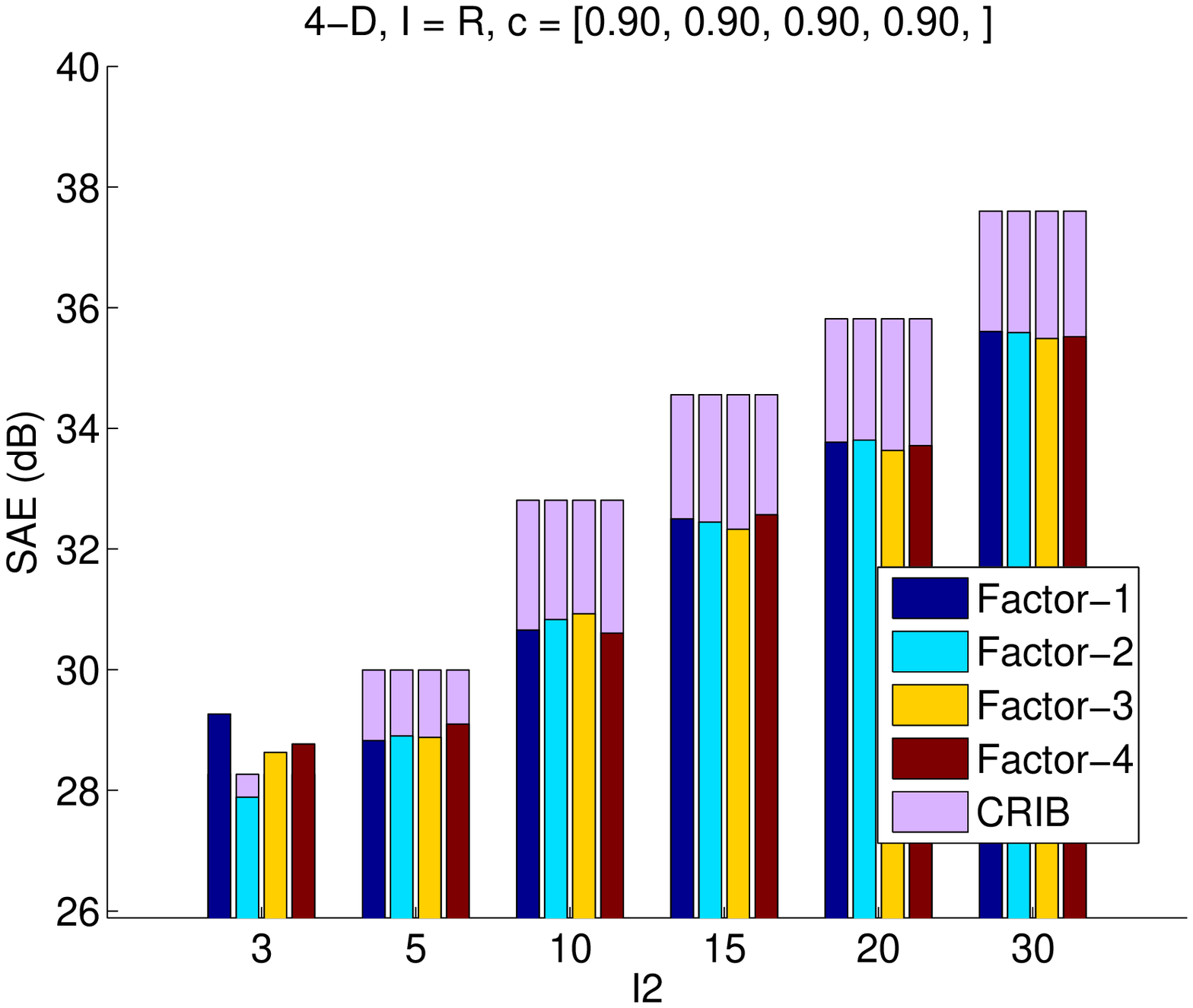}\label{fig_N4_c9999_msaevscrib_allcomps}}
\vspace{-1ex}
\caption{Median SAE of all components for factors over 30 Monte Carlo runs vs CRIB in decomposition of order-4 tensors of size $I_n = R$ = 3, 5, 10, \ldots, 30, for all $n$ through the unfolding rule $\bl = [1, 2, (3,4)]$.
Correlation coefficients $c_n$ have been chosen from the set $ \in \{0.1, 0.9\}$, for all $n$.
The signal to white Gaussian noise power ratio (SNR) was at 10dB or 30dB.
}\label{fig_4D_lossfactors}
\end{figure}

\begin{example}\label{ex_1}
\end{example}
We illustrate the similar behavior of CRIB over unfolding but for higher-order ranks.
We decomposed $\tY_{\llbracket 1,2,(3,4) \rrbracket}$ unfolded from rank-$R$ tensors of size $R \times R \times R \times R$ with $R  =  3, 5, \ldots, 30$, corrupted with additive Gaussian noise of 10 dB SNR. 
There was not any significant loss in factors when modes 1 and 2 comprised low-collinear components despite of collinearity in modes 3 and 4 as seen in Figs.~\ref{fig_N4_c1111_msaevscrib_allcomps}-\ref{fig_N4_c1199_msaevscrib_allcomps}.
For all the other cases of $(c_1,c_2,c_3,c_4)$, there were always significant losses, especially when all the factors comprised highly collinear components (i.e., $c_n$ is close to $\pm1$) as seen in Figs.~\ref{fig_N4_c1919_msaevscrib_allcomps}-\ref{fig_N4_c9999_msaevscrib_allcomps}.

When the first mode has nearly collinear factors, i.e.,  $c_1$ is close to $\pm 1$, we have \cite{2012arXiv1209.3215T}
\begin{eqnarray}
\mbox{CRIB}(\ba_1)_{c_1=\pm 1}&=&\frac{\theta(I_1-1)}{1-h^2}
<\mbox{CRIB}(\ba_1)_{c_1=0},
\end{eqnarray}
but the expressions for the folded tensor decomposition remain
unchanged. It means that the loss occurs as seen in Fig.~\ref{fig_N4_c1919_msaevscrib_allcomps} and Fig.~\ref{fig_N4_c9999_msaevscrib_allcomps}.

\subsubsection{Loss in unfolding order-5 tensors}\label{sec::loss5D}
 For order-5 rank-2 tensors, we consider the case when $c_1=0$, and put $h=c_2c_3c_4c_5$. CRIBs of decompositions of the full and unfolded tensors are given by
\be
\mbox{CRIB}(\ba_1)&=&\frac{\theta}{1-h^2}\left[I_1-1+\frac{\zeta-4h^2}
{1+3h^2-\zeta}\right] , \label{(7)}\\
%
\mbox{CRIB}_{\llbracket 1,2,(3,4,5) \rrbracket}(\ba_1)&=&\frac{\theta}{1-h^2}\left[I_1-3+\frac{1}{1- c_2^2} + \frac{1}{1-c_3^2c_4^2c_5^2}\right], \label{(8)}\\
\mbox{CRIB}_{\llbracket 1,(2,3),(4,5) \rrbracket}(\ba_1)&=&\frac{\theta}{1-h^2}\left[I_1-3+\frac{1}{1-c_2^2c_3^2} + \frac{1}{1-c_4^2c_5^2}\right], \label{(9)}
\ee
where $\zeta = c_2^2c_3^2c_4^2+c_2^2c_3^2c_5^2+c_2^2c_4^2c_5^2+c_3^2c_4^2c_5^2 $. 
From (\ref{(8)}) and (\ref{(9)}), it is obvious that $\mbox{CRIB}_{\llbracket 1,2,(3,4,5) \rrbracket}(\ba_1) \le \mbox{CRIB}_{\llbracket 1,(2,3),(4,5) \rrbracket}(\ba_1)$  if $c_2^2 \le c_4^2c_5^2$. 
This rule coincides with that for order-4 tensors to reduce the collinearity of the merging factor matrices as much as possible.
For $c_2=0$, the expressions (\ref{(7)}) and
(\ref{(8)}) become identical, but expression (\ref{(9)}) is
larger, in general. 

\begin{figure}
\centering
\psfrag{c}[t][t]{\scalebox{1}{\color[rgb]{0,0,0}\setlength{\tabcolsep}{0pt}\begin{tabular}{c}\\[-2.5em]\small $c$\end{tabular}}}%
\psfrag{Correlation coefficient}[t][t]{\scalebox{1}{\color[rgb]{0,0,0}\setlength{\tabcolsep}{0pt}\begin{tabular}{c}\\[-2.5em]\small $c$\end{tabular}}}%
\subfigure[CRIB (dB) for order-6 rank-2 tensors.]
{
\psfrag{CRIB}[l][l]{\scalebox{.8}{\color[rgb]{0,0,0}\setlength{\tabcolsep}{0pt}\vspace{3ex}\begin{tabular}{c}\footnotesize $\mbox{CRIB}$\end{tabular}}}%
\psfrag{CRIBl1}[l][l]{\scalebox{.8}{\color[rgb]{0,0,0}\setlength{\tabcolsep}{0pt}\vspace{3ex}\begin{tabular}{c}\footnotesize $\mbox{CRIB}_{\bl_1}$\end{tabular}}}%
\psfrag{CRIBl2}[l][l]{\scalebox{.8}{\color[rgb]{0,0,0}\setlength{\tabcolsep}{0pt}\vspace{3ex}\begin{tabular}{c}\footnotesize $\mbox{CRIB}_{\bl_2}$\end{tabular}}}%
\psfrag{CRIBl3}[l][l]{\scalebox{.8}{\color[rgb]{0,0,0}\setlength{\tabcolsep}{0pt}\vspace{3ex}\begin{tabular}{c}\footnotesize $\mbox{CRIB}_{\bl_3}$\end{tabular}}}%
\psfrag{CRIBl4}[l][l]{\scalebox{.8}{\color[rgb]{0,0,0}\setlength{\tabcolsep}{0pt}\vspace{3ex}\begin{tabular}{c}\footnotesize $\mbox{CRIB}_{\bl_4}$\end{tabular}}}%
\psfrag{CRIBl5}[l][l]{\scalebox{.8}{\color[rgb]{0,0,0}\setlength{\tabcolsep}{0pt}\vspace{3ex}\begin{tabular}{c}\footnotesize $\mbox{CRIB}_{\bl_5}$\end{tabular}}}%
\includegraphics[width=.48\linewidth, trim = 0.0cm 0cm 0cm 0cm,clip=true]{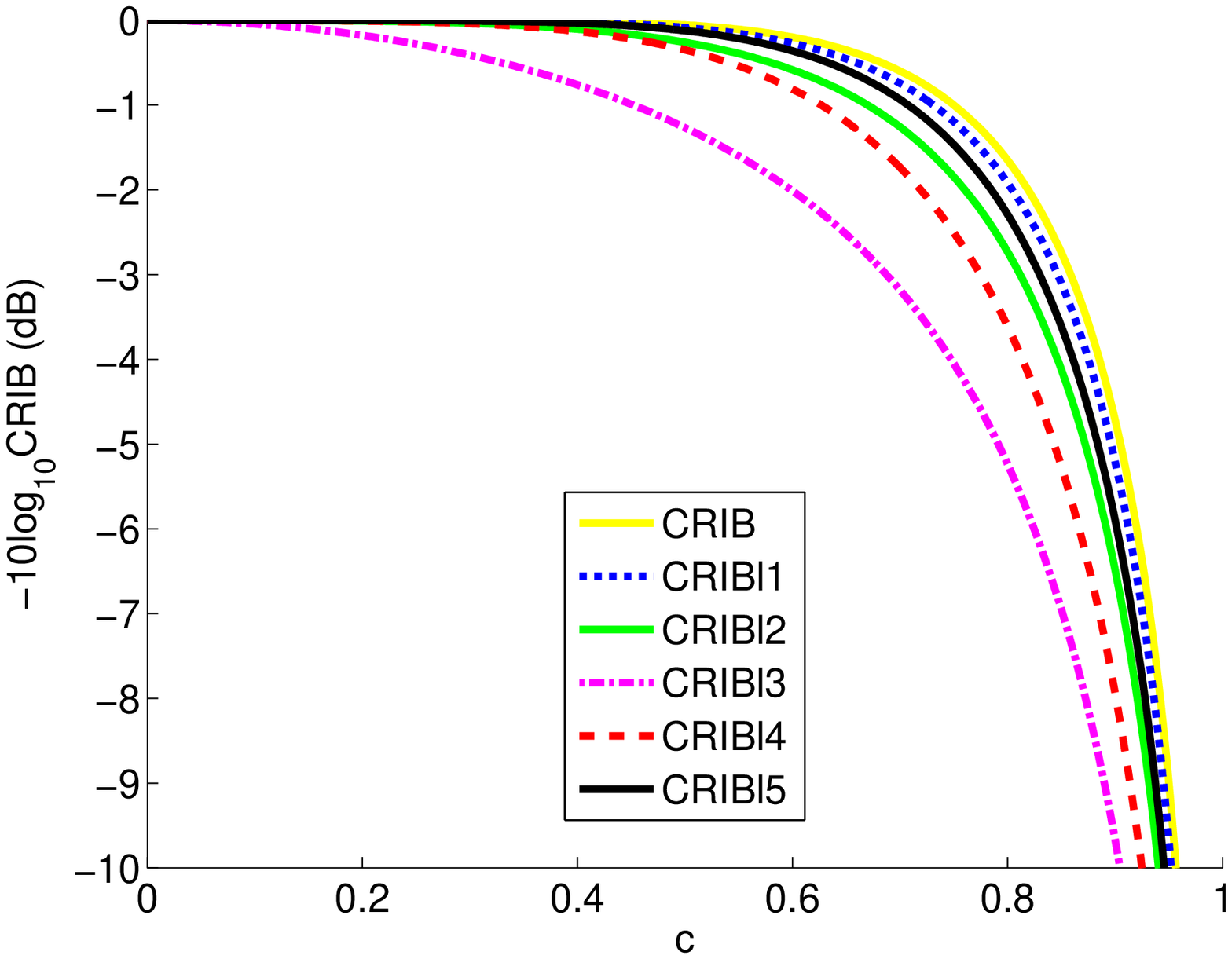}\label{fig_rnk2_ord6}}
\hfill
%
%
%
%
%
\subfigure[CRIB loss for rank-20 tensors.]
{
\includegraphics[width=.47\linewidth, trim = 0.0cm 0cm 0cm 0cm,clip=false]{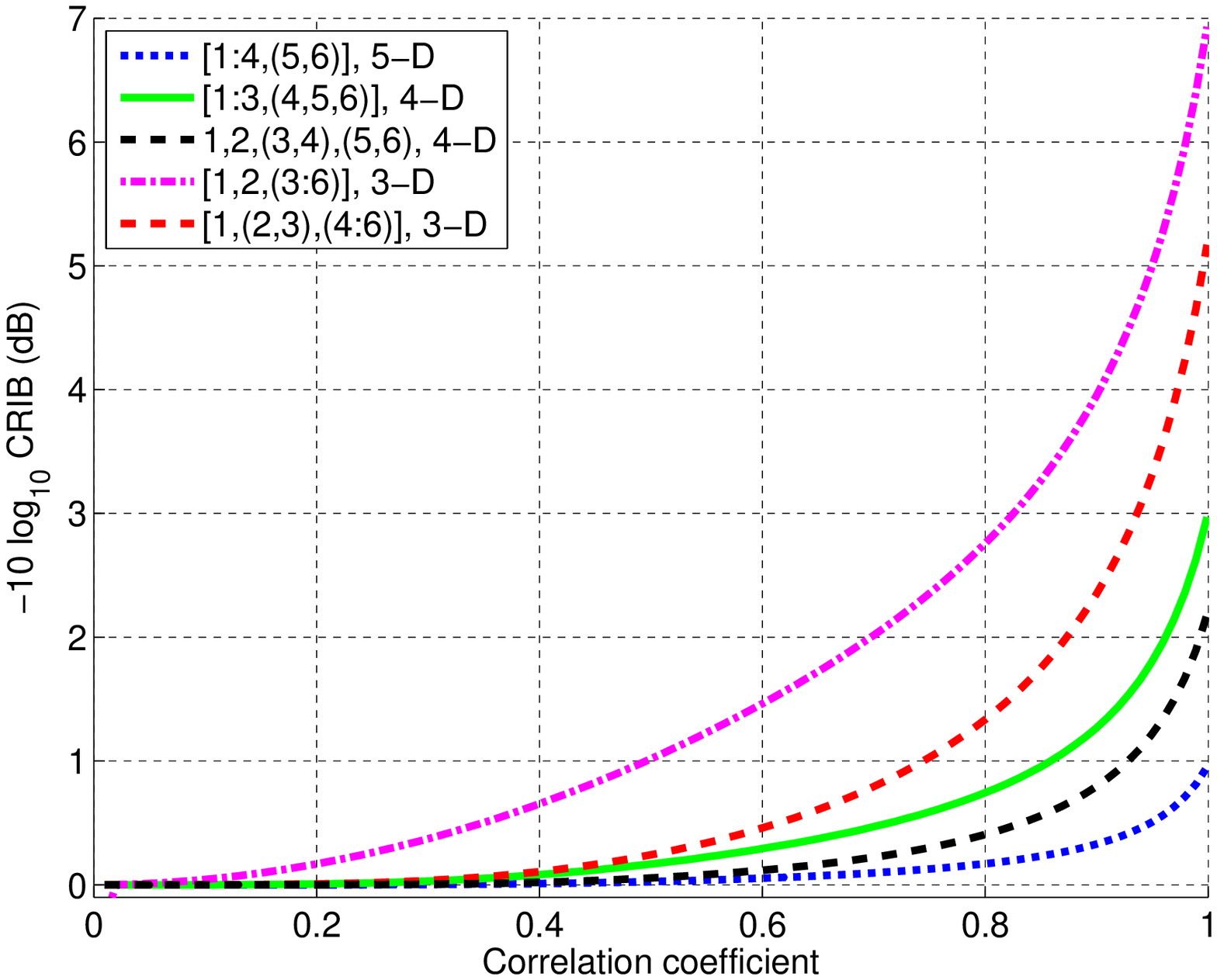}\label{fig_loss_crib_folding_R20}}
%
%
\caption{
%
%
The CRIB loss in decomposition of order-6 rank-$R$ tensors with size $I_n = R$ and correlation coefficients $c_n = c$, for all $n$ following 5 unfolding rules.
The CRIB loss is significant when components are collinear, i.e., $c \rightarrow 1$. Unfolding $\bl = [1{\emph{:}}4,(5,6)]$ causes a lesser CRIB loss than other rules. The unfolding $\bl = [1,2,(3,4),(5,6)]$ is more efficient than multimode unfolding rule.
}
\label{fig_CRIB_N6R5_fold3}
\end{figure}

For higher order tensors, we analyze the CRIB loss in decomposition of order-6 tensors with assumption that $c_1 = c_2 = \cdots = c_6 = c$ through 5 unfoldings $\bl_1 = [1,2,3,4,(5,6)]$, $\bl_2 = [1,2,3,(4,5,6)]$, $ \bl_3 = [1,2,(3,4,5,6)]$, $\bl_4 = [1,(2,3),(4,5,6)]$ and $\bl_5 = [1,2,(3,4),(5,6)]$
\be
\mbox{CRIB}(\ba_1) &=& \theta\left[\frac{I_1-1}{1-c^{10}} + \frac{5c^8(4c^6+3c^4+2c^2+1)}{(1-c^{10})(1-c^8)(1+3c^2+c^4)(1+c^2+6c^4+c^6+c^8)}\right] ,  \notag \\
\mbox{CRIB}_{\bl_1}(\ba_1) &=& \theta\left[\frac{I_1-1}{1-c^{10}} + \frac{c^6(6c^8+11c^6+7c^4+5c^2+1)}{(1-c^{10})(1-c^8)(1+c^2)(1+2c^2+6c^4+2c^6+c^8)}\right] ,  \notag \\
\mbox{CRIB}_{\bl_2}(\ba_1) &=& \theta\left[\frac{I_1-1}{1-c^{10}} + \frac{c^4(4c^6+3c^4+2c^2+1)}{(1-c^{10})(1-c^8)(1+3c^2+c^4)}\right] ,  \notag \\
\mbox{CRIB}_{\bl_3}(\ba_1) &=& \theta \left[\frac{I_1-1}{1-c^{10}} + \frac{c^2(2c^6+c^4+c^2+1)}{(1-c^{10})(1-c^8)}\right] ,  \notag \\
 \mbox{CRIB}_{\bl_4}(\ba_1) &=& \theta \left[\frac{I_1-1}{1-c^{10}} + \frac{c^4(1+c^4)(2c^4+2c^2+1)}{(1-c^{10})(1-c^8)(1+c^2+c^4)(1-c+c^2)}\right] ,  \notag \\
 \mbox{CRIB}_{\bl_5}(\ba_1) &=& \theta \left[\frac{I_1-1}{1-c^{10}} + \frac{c^6(2c^6 + 2c^4 + 2c^2 + 1)(c^6 + 4c^4 + 3c^2 + 2)}{(1-c^{10})(1-c^8)(1+c^2+c^4)(c^8+3c^6+6c^4+3c^2+1)}\right] .  \notag
 \ee
It holds $\mbox{CRIB} <  \mbox{CRIB}_{\bl_1} < \mbox{CRIB}_{\bl_5} < \mbox{CRIB}_{\bl_2} < \mbox{CRIB}_{\bl_4} < \mbox{CRIB}_{\bl_3}$ (see Fig.~\ref{fig_rnk2_ord6} for $\theta = 1$ and $I_1 = R$). 
%
%
%
%
%
Fig.~\ref{fig_loss_crib_folding_R20} illustrates such behavior of the CRIB losses $\calL = -10\log_{10}\left(\frac{\mbox{CRIB}(\ba_1)}{{\mbox{CRIB}_{\bl}(\ba_1)}}\right)$ (dB) but for higher rank $R = 20$ and the tensor size $I_n = R = 20$, for all $n$. 
The CRIB loss is insignificant when components are nearly orthogonal (i.e., $c \rightarrow 0$), and is relevant for highly collinear components, i.e., $c \rightarrow 1$. 
Unfolding $\bl_1$ causes a CRIB loss less than 1 dB, while unfolding $\bl_2$, $\bl_4$ and $\bl_3$ can cause a loss of 3, 5 and 7 dB, respectively.
The result confirms that two-mode unfolding causes a lesser CRIB loss than other rules.
The unfoldings $\bl_4$ and $\bl_5$ are more efficient than multimode unfoldings $\bl_2$ and $\bl_3$, respectively in decomposition of unfolded tensors of the same orders.


\subsubsection{A case when two factor matrices have orthogonal columns}\label{sec::orthoCP}

As pointed out in (\ref{eq_ortho4}) that there is not any loss of accuracy in estimation of $\bA^{(1)}$ and $\bA^{(2)}$ through unfolding when these two factor matrices have mutually orthogonal columns. 
The result also holds for arbitrary order-$N$ rank-$R$ tensors which have orthogonal components on two modes.
In such case, the analytical CRIB is given by
\begin{theorem}[\cite{2012arXiv1209.3215T}] \label{theo_orthoCP} When $\bA^{(1)}$ and $\bA^{(2)}$ have mutually orthogonal columns, it holds
\be
\mbox{CRIB}(\ba_1)&=&  \frac{\sigma^2}{\lambda_1^2} \left(I_1-R+\sum_{r = 2}^{R}\frac{1}
{1-\gamma_{r}^2}\right),\; \quad \gamma_{r}  = \prod_{n = 3}^{N} \left(\ba^{(n)\, T}_1 \ba^{(n)}_r\right), \; r = 2, 3, \ldots, R. \label{eq_criborthoN}
\ee
\end{theorem}
It is obvious that $\mbox{CRIB}_{\llbracket1,2,(3\emph{:}N)\rrbracket}(\ba_1) = \mbox{CRIB}(\ba_1)$.
Hence, estimation of $\bA^{(1)}$ and $\bA^{(2)}$ through unfolding is lossless in terms of accuracy.

An important observation from Theorem~\ref{theo_orthoCP} is that  
all the factor matrices in CPD with orthogonality constraints\cite{Dosse_indort,Sorensen_orthoCPD} can be estimated through order-3 tensors unfolded from the original data without any loss of accuracy. That is an algorithm for order-3 orthogonally constrained CPD on two or three modes can work for arbitrary order-$N$ tensors.

\subsection{Unfolding strategy}\label{sec::unfoldingstrategy}
Based on the above analysis of CRIB loss, we summarize and suggest an efficient unfolding strategy to reduce the loss of accuracy.
Without loss of generality, assume that $0 \le |c_1| \le |c_2| \le \cdots \le  |c_N| \le 1$, the following procedures should be carried out when unfolding an order-$N$ tensor to order-$M$ (typically, $M$ = 3)
\begin{itemize}
\item Unfold two modes which correspond to the two largest values $|c_n|$, i.e., $(N-1,N)$. This yields a new factor matrix with a correlation coefficient $\tilde{c}_{N-1} = c_{N-1} c_{N}$. The tensor order is reduced by one.
\item Unfold two modes which correspond to the two largest collinearity values among $(N-1)$ values $[c_1, c_2, \ldots, c_{N-2}, \tilde{c}_{N-1}]$. 
This can be $(N-3,N-2)$ if $|\tilde{c}_{N-1}| < |c_{N_2}|$; otherwise,  $(N-2,N-1,N)$.
The new correlation coefficient is $c_{N-3} c_{N-2}$ or $c_{N-2} c_{N-1} c_{N}$.
\item  Unfold the tensor until its order is $M$.
\end{itemize}
In addition, (nearly) orthogonal modes should not be merged in the same group.
For order-4 tensor, the unfolding [1,2,(3,4)] is recommended. 

\begin{figure}
\psfrag{I2}[t][b]{\scalebox{1}{\color[rgb]{0,0,0}\setlength{\tabcolsep}{0pt}\vspace{3ex}\begin{tabular}{c}\small $I$\end{tabular}}}%
\centering
{
\includegraphics[width=.48\linewidth, trim = 0.0cm -.5cm 0cm 1cm,clip=true]{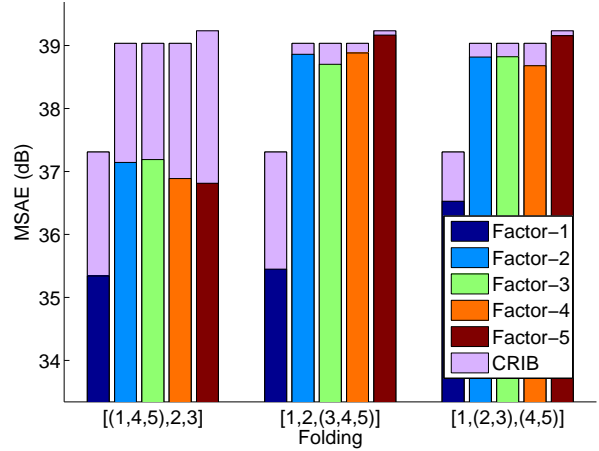}\label{fig_N3_c117778_cribfoldvscrib_allcomps}}
\caption{Affect of unfolding rules to the accuracy loss in decomposition of order-5 tensors of size $I_n = R = 10$ with $c_1 = 0.1$, $c_2 = c_3 = c_4 = 0.7$ and $c_5 = 0.8$. Mean SAEs (dB) were computed for all the components over 100 Monte Carlo runs.}
\end{figure}

\begin{example}\label{ex_3}
\end{example} We decomposed order-5 tensors with size $I_n = R = 10$ and additive Gaussian noise of 10 dB SNR. Correlation coefficients of factors matrices were [0.1, 0.7, 0.7, 0.7, 0.8]. 
Three tensor unfoldings $\bl_1 = [(1,4, 5), 2, 3]$, $\bl_2 = [1, 2, (3, 4, 5)]$ and $\bl_3 = [1, (2,3), (4,5)]$ were applied to the order-5 tensors. 
Unfolding $\bl_1 =  [(1,4, 5), 2, 3]$ caused the largest loss of 2 dB with an average MSAE = 36.62 dB illustrated in Fig.~\ref{fig_N3_c117778_cribfoldvscrib_allcomps}. The  recommended unfolding $\bl_3$ according to the above strategy achieved an average MSAE = 38.29 dB compared with the average CRIB = 38.67 dB on all the components.

\begin{example}{\textbf{Unfolding tensors with the same colinearity in all modes.}}\label{ex_4}
\end{example}  We verified the unfolding rules for order-6 tensors with simplified assumption that $c_1 = c_2 = \cdots = c_6 = c$.
Since correlation coefficients are identical, the unfoldings $\bl = [1, 2, (3, 4), (5,6)]$ is one of the best rules.
Fig.~\ref{fig_loss_crib_folding_R20} shows that the CRIB loss by $\bl = [1, 2, 3, (4, 5, 6)]$ was higher than the loss by $\bl = [1, 2, (3, 4), (5,6)]$.

\subsection{Unfolding without collinearity information}\label{sec::unfoldingrealdata}

For real-world data, although collinearity degrees of factor matrices are unknown, the above strategy is still applicable. 
Since the decomposition through tensor unfolding decomposes only an order-3 tensor, the computation is relatively fast. We can try any tensor unfolding, and verify the (average) collinearity degrees of the estimated factors $\displaystyle c_n = \frac{\sum_{r \neq s} |\ba^{(n) T}_r \ba^{(n)}_s|}{R(R-1)}$, ($n = 1, 2, \ldots, N$), to proceed a further decomposition with a new tensor unfolding.

%
%

%
%
%
%
%

{\vspace{10pt}\par\noindent{\bf Example \ref{ex_3}(b)}} We replicated simulations in Example~\ref{ex_3}, but assumed that there was not prior collinearity information and the bad unfolding rule $\bl_1 =  [(1,4, 5), 2, 3]$ was applied. The average collinear degrees of the estimated factor matrices $c_n$ = 0.0989,  0.7007,  0.6992,   0.7021, 0.8014, for $n = 1, \ldots, N$, respectively,  indicated that the unfolding $[(1,4, 5), 2, 3]$ is not a good one.
The unfolding $\bl_3 = [1, (2,3), (4,5)]$ was then suggested, and a further decomposition was proceeded. This improved the MSAE about 2dB.  

For more examples, see decomposition of the ITPC tensor in Example~\ref{ex_8} when $R = 8$.


\section{Fast Approximation for High Order and Difficult scenario CPD}\label{sec::3wtoNw-lowrank}



An appropriate unfolding rule can reduce the loss of the decomposition. However, the loss always occurs when factors have high collinearity or unfolding orthogonal modes.
Moreover, in practice, a rank-$R$ CPD may not fit the data tensor.
This could happen when $R$ is not exactly the true rank of the tensor.
Especially, for under-rank CPD, the error tensor $\tE$ can still be explained by tensors in the Kruskal form. 
In this case, components of the merging factor matrices tend to 
comprise information of the other components in higher rank CP decomposition. Hence, they are no longer rank-one matrices/tensors, and approximation of merging components by rank-one tensors cannot yield good approximate to the true factors.
To this end, low-rank approximations to merging components are proposed, and  components are estimated through two major stages
\begin{enumerate}
\item Construct an order-$N$ structured Kruskal tensor $\tensor{\tilde Y}_{J}$ from the order-$M$ rank-$R$ Kruskal tensor which approximates the unfolded tensor $\tY_{\llbracket \bl \rrbracket}$. $\tensor{\tilde Y}_{J}$ often has higher rank than $R$.
\item Approximate $\tensor{\tilde Y}_{J}$ by a rank-$R$ Kruskal tensor which is the  final result.
\end{enumerate}

The algorithm is first derived for unfolding two modes, and extended to multimode unfolding.



\subsection{Unfolding two modes}\label{sec::fold2modes}


We consider a rank-$R$ CPD of $\tY$ and a simple unfolding $\bl = [1, \ldots,N-2, (N-1, N)]$ 
\be
	\tY_{\llbracket  \bl \rrbracket} = \sum_{r = 1}^{R} \lambda_r \, \bb^{(1)}_r  \circ \bb^{(2)}_r \circ \cdots \circ \bb^{(N-1)}_r + \tE. \label{equ_folding_rnkR_4D}
\ee
Assume matrices $\bF_{r} = {\tt{reshape}}(\bb^{(N-1)}_r,{[I_{N-1} \times I_N ]})$ have rank-$J_r$ ( $1\le J_r \ll I_{N-1}$), i.e., 
$\bF_{r} = \bU_r  \bSigma_r \bV_r^T$, for $r = 1, 2, \ldots, R$,
where $\bSigma_r = \diag(\bsigma_r)$ and singular values $\bsigma_r = [\sigma_{r1}, \sigma_{r2},\ldots, \sigma_{rJ_r}]$, $1\ge \sigma_{r1} \ge \sigma_{r2} \ge \cdots \ge \sigma_{rJ_r} > 0$, $\sum_{j= 1}^{J_r} \sigma_{rj}^2  = 1$. By replacing all $\bb^{(N-1)}_r$ by matrices $\bU_r$ and $\bV_r$ for $r = 1,2, \ldots, R$, and replicating components $\bb^{(n)}_r$ ($n = 1, 2, \ldots, N-1$) $J_r$ times, we generate an order-$N$ rank-$J$ tensor in the Kruskal form where $J = \sum_r J_r$.


\begin{lemma}\label{lem_rankJ_CP}
The order-$N$ rank-$J$ Kruskal tensor $\tensor{\widetilde Y}_J = \llbracket  {\tilde{\blambda}}; \tilde\bA^{(1)} , \ldots,  \tilde\bA^{(N-1)} ,   \tilde\bA^{(N)}\rrbracket$,  where $\displaystyle J = \sum_{r = 1}^{R} J_r$, ($R\le J \ll R \min(I_{N-1},I_N)$, 
$\tilde\blambda = 	 \left[ 
	 \begin{array}{*{10}c}
	 	\lambda_1  \1_{J_1}^T& \lambda_2 \1_{J_2}^T & \ldots & \lambda_R \1_{J_R}^T
	 \end{array} 
\right]^T   \;  \in \Real^{J}$, and 
\be
\tilde\bA^{(n)} &=& \begin{cases}
	 \bB^{(n)}   \,  \bM \quad  \in \Real^{R \times J}, \quad & n = 1, \ldots, N-2, \notag\\
	 \left[ 
	 \begin{array}{*{10}c}
	 	\bU_1 & \bU_2 & \ldots & \bU_R
	 \end{array} 
	\right]   \quad  \in \Real^{I_3 \times J} ,  & n = N-1, \notag\\
	 \left[ 
	 \begin{array}{*{10}c}
	 	\widetilde\bV_1 & \widetilde \bV_2 & \ldots & \widetilde\bV_R  
	 \end{array} 
	\right]   \quad  \in \Real^{I_4 \times J} , \quad \widetilde\bV_r = \bV_r \diag(\bsigma_r)  & n = N,
	\end{cases} \\
	\bM &=& \blkdiag(\1_{1\times J_1}, \1_{1\times J_2}, \ldots, \1_{1\times J_R}),\notag
\ee
has the same approximation error of the best rank-$R$ CPD of $\tY_{\llbracket \bl \rrbracket}$, i.e.,
$
	\|\tY - \tensor{\widetilde Y}_J\|_F^2 =  \|\tE\|_F^2.
$
\end{lemma}

If $J_r = 1$ for all $r$, $\tensor{\widetilde Y}_J$ is the approximation ${\tensor{\widehat{Y}}}$ of the true factors as pointed out in the previous section.
Otherwise, the order-$N$ rank-$J$ Kruskal tensor $\tensor{\widetilde Y}_J$ is approximated by a rank-$R$ Kruskal tensor ${\tensor{\widehat{Y}}}$.
Note that this procedure does not access and manipulate on the real data $\tensor{\widetilde Y}$.
For example, the mode-$n$ CP-gradients of $\tensor{\widetilde Y}_J$ with respect to $\bA^{(k)}$ ($k\neq n$, $n = 1, 2, \ldots, N$) which are the largest workload in CP algorithms such as ALS, OPT, dGN can be quickly computed as illustrated in Appendix~\ref{sec::app_fastCP}
with a low computational complexity of order $\displaystyle \mathcal O\left(JR(I_{N-1} + I_{N}) +  R^2 \left(\sum_{n= 1}^{N-2} I_n + N\right)\right)$.
It means that computation of $\left(\widetilde\bY_{J(n)} - \hat\bY_{(n)}\right) \, \left(\bigodot_{k\neq n} \bA^{(k)} \right)$ is much faster than the computation on the raw data
$\left(\bY_{(n)} - \hat\bY_{(n)}\right) \, \left(\bigodot_{k\neq n} \bA^{(k)} \right)$. In other words, estimation of factors $\bA^{(n)}$ from the Kruskal tensor $\tensor{\widetilde Y}_J$ is relatively fast.
%

When the matrices $\bF_{r}$ ($r = 1, 2, \ldots, R$) have not exact low-rank representation, we consider their truncated SVDs
such that 
$\displaystyle \rho_r =   \sum_{j= 1}^{J_r} \sigma_{rj}^2  \ge  \tau$ , $0 \ll \tau \le 1$.
The parameter $\tau$ controls the loss of fit by low rank approximations. 
The higher $\tau$, the lower loss of fit, but the higher the approximation ranks.
In the simulations, we set $\tau \ge 0.98$.

Let $\tensor{\widetilde Y}_{R}$ denote the solution of the rank-one FCP algorithm (i.e., using Algorithm~\ref{alg_3w2Nw_rank1})
\be
	\tensor{\widetilde Y}_{R} &=&  \left\llbracket \blambda_{R}; \bB^{(1)}, \ldots,  \bB^{(N-2)}, \left[\bu_{11}, \bu_{21}, \ldots, \bu_{R1}\right],  \left[\bv_{11}, \bv_{21}, \ldots, \bv_{R1}\right] \right\rrbracket, 
\ee
where $\blambda_{R} = [\lambda_1 \sigma_{11}, \lambda_2 \sigma_{21}, \ldots, \lambda_R \sigma_{R1}] \in \Real^{R}$.
It is straightforward to see that
$
	\bb^{(n-1)}_r = \sum_{q = 1}^{\min(I_{N-1},I_N)} \sigma_{rq} \, ( \bv_{rq} \otimes \bu_{rq})$, $r = 1, 2, \ldots, R$.
Because $\bB^{(n)} = [\bb^{(n)}_1, \ldots, \bb^{(n)}_R]$, ($n = 1, \ldots, N-1$), forms the best rank-$R$ CPD of $\tY_{\llbracket \bl \rrbracket}$, each vector $\sigma_{rq} (\bv_{rq} \otimes \bu_{rq})$, ($r = 1, \ldots, R$), contributes to achieve the optimal approximation error $\|\tE\|_F$ in (\ref{equ_folding_rnkR_4D}).
Discarding any set of singular components $(\bv_{rq}  \otimes \bu_{rq})$ will increase the approximation error. 
The more singular vectors to be eliminated, the higher approximation error of $\tY_{\llbracket \bl \rrbracket}$.
%
It means that the tensor $\tensor{\widetilde Y}_{R}$ has higher approximation error than 
the tensor $\tensor{\widetilde Y}_{J}$. That is
\be
	\|\tE\|_F^2 \le \|\tY - \tensor{\widetilde Y}_J\|_F^2 \le \|\tY - \tensor{\widetilde Y}_R\|_F^2, \label{equ_bounderror_folding2}
\ee
or performance of FCP using low-rank approximations is  better than that using Algorithm~\ref{alg_3w2Nw_rank1}.

\subsection{Unfolding $M$ modes}\label{sec::foldMmode}

Consider an example where unfolding $M$ modes are $\bl = [1, \ldots, N-M, (N-M+1, \ldots, N)]$.
In this case, truncated SVD is no longer applied to approximate tensors $\tF_{r}$ of size ${I_{N-M+1} \times \cdots \times I_N }$ and $\vtr{\tF_{r}} = \bb^{(N-M)}_r$.
However, we can apply low-rank Tucker decomposition to $\tF_{r}$
\be
	\tF_{r} \approx \llbracket \tS_r  ;  \bU^{(1)}_{r},  \bU^{(2)}_{r},  \ldots ,  \bU^{(M)}_{r} \rrbracket, \quad (r = 1, \ldots, R),\label{equ_MSVD_Br}
\ee
where $\bU^{(m)}_{r}$ are orthonormal matrices of size $(I_{N-M+m-1} \times {T_{rm}})$, and the core tensor $\tS_r = [s^{(r)}_{\bt}]$ is of size ${T_{r1}} \times {T_{r2}} \times \cdots \times {T_{rM}}$, $\bt = [t_1, t_2, \ldots, t_M]$.

In order to estimate an order-$N$ rank-$R$ Kruskal tensor $\tensor{\widehat Y}$, Tucker tensors in (\ref{equ_MSVD_Br}) are converted to an equivalent Kruskal tensor of rank-$(T_1T_2 \ldots T_M)$. 
However, we select only the most $J_r$ dominant $s^{(r)}_{\bt}$ ($1\le J_r \ll T_1T_2 \ldots T_M$), ${\bt} \in \mathcal T = \{{\bt}_1, {\bt}_2,\ldots, {\bt_{J_r}}\}$ 
among all coefficients of $\tS_r$  so that
\be
	\rho_r = \sum_{\bt  \in \mathcal T} (s^{(r)}_{\bt}  )^2  \ge \tau \, \|\tF_r\|_F^2.
\ee
The tensors $\tF_{r}$ have rank-$J_r$ approximations in the Kruskal form
\be
	\tB_{r} \approx \llbracket \blambda_r;  \bU_{r}^{(1)}   \bM_{r1} , \bU_{r}^{(2)} \,\bM_{r2},  \ldots,  \bU_{r}^{(M)}  \bM_{rM}\rrbracket, \quad (r = 1, 2,\ldots, R),\notag
\ee
where $\blambda_r  = [s^{(r)}_{\bt_1}, s^{(r)}_{\bt_2}, \ldots, s^{(r)}_{\bt_{J_r}}]$ and $\bM_{rm}$ ($m = 1, 2, \ldots, M$) are indicator matrices of size $T_m \times J_{r}$ which have only $J_r$ non-zero entries at $\bM_{rm}(t_{j,m}, j) = 1$, $\bt_{j} = [t_{j,1}, t_{j,2}, \ldots, t_{j,M}] \in \mathcal T$, $j = 1, \ldots, J_r$.

Combination of $R$ rank-$J_r$ CP approximations for components $\bb^{(N-1)}_r$ yields a rank-$J$ Kruskal tensor $\tensor{\widetilde Y}_J$ ($J = \sum_{r}^{R} J_r$) as mentioned in the previous section (see Lemma~\ref{lem_rankJ_CP}). 
A rank-$R$ CPD of $\tensor{\widetilde Y}_J$ will give us an approximate to the true solution.

An alternative approach is that we consider $M$-mode unfolding as $(M-1)$ two mode unfoldings. For example, since $(1,2,3) \equiv (1,(2,3))$, the factor matrices are then sequentially estimated using the method in Section \ref{sec::fold2modes}. 
Indeed, this sequential method is recommended because it is based on SVD and especially low-rank approximation to matrix is well-defined.

\def\baselinestretch{1} 
\setlength{\algomargin}{1em}
\begin{algorithm}[ht]
\SetFillComment
\SetSideCommentRight
\caption{{\tt{FCP}}\label{alg_3w2Nw_lowrank}}
\DontPrintSemicolon \SetFillComment \SetSideCommentRight
\KwIn{Data tensor $\tY$:  $(I_1 \times I_2 \times \cdots \times I_N)$, rank $R$, threshold $\tau (\ge 0.98)$\linebreak
	Unfolding rule $\bl = [\bl_1, \bl_2, \ldots, \bl_M]$, $\bl_m = [l_m(1), l_m(2), \ldots, l_m(K_m)]$} 
\KwOut{$\blambda \in \Real^{N}$, $N$ matrices $\bA^{(n)} \in \Real^{I_n \times R}$} \SetKwFunction{mreshape}{reshape}
\SetKwFunction{permute}{permute} 
\SetKwFunction{svds}{svds} 
\SetKwFunction{eigs}{eigs} 
\SetKwFunction{rankone}{rank1CP} 
\SetKwFunction{lowrank}{lowrankCP} 
\SetKwFunction{CPD}{CPD} 
\SetKwFunction{psCPD}{structuredCPD} 
\SetKwFunction{MwtoNw}{FCP} 
\SetKwFunction{tucker}{TD} 
\Begin{
 {\mtcc{Stage 1: Tensor unfolding and compression \dashrule}}
\nl $ \llbracket \tG, \bU^{(1)},  \ldots, \bU^{(M)} \rrbracket = \tucker(\tY_{\llbracket \bl \rrbracket},\min(\bI, R))${\mtcc*{Tucker decomposition of order-$M$ $\tY_{\llbracket \bl \rrbracket}$}}
\;\vspace{-3ex}
{\mtcc{Stage 2: CPD of the unfolded tensor \dashrule}}
{
\nl $ \llbracket \blambda; \bB^{(1)},  \ldots, \bB^{(M)} \rrbracket = \CPD(\tG,R)${\mtcc*{order-$M$ CPD of the core tensor }}
}
\nl \lFor{$m  = 1, 2, \ldots, M$}{
$\bA^{(m)} \leftarrow \bU^{(m)} \, \bB^{(m)}${\mtcc*{Back projection of TD}}
}
\;\vspace{-3ex}
{\mtcc{Stage 3: Sequential estimation of factor matrices  \dashrule}}
\nl $\widetilde M = M$\;
  \For{$m  = M, M-1, \ldots, 1$}{
  \If{$K_m \ge 2$}{
  \For{$k  = 1, 2, \ldots, K_m$}{
   {\mtcc{Stage 3a: Construction of rank-$J$ Kruskal tensor \dashrule}}
   
\nl ${\widetilde m} = m + k-1$, $\tilde\blambda = []$, $\tilde\bA^{(\widetilde m)} = []$, $\breve\bA^{(\widetilde m)} = []$, $\bM = []$\;
 \For{$r  = 1, 2, \ldots, R$}{
\nl $\bF_{r} = {\tt{reshape}}(\ba^{({\widetilde m})}_r,{[I_{l_m(k)}, \prod_{i = k+1}^{K_m} I_{l_m(i)}]})$\;
\nl $\bF_{r}  \approx \bU_r \diag(\bsigma_r) \, \bV_r^T${\mtcc*{truncated SVD such that $\|\bsigma_r\|_2^2 \ge \tau$}}
\nl   $\tilde\blambda  \leftarrow [ \tilde\blambda , \lambda_r \, \bsigma_r]$,\quad
   $\tilde\bA^{({\widetilde m})}  \leftarrow [\tilde\bA^{({\widetilde m})},   \bU_r]$,\quad
   $\tilde\bA^{({\widetilde m}+1)}  \leftarrow [\tilde\bA^{({\widetilde m}+1)},   \bV_r]$\;
 \nl   $\blambda_r  \leftarrow  \lambda_r \, \sigma_{r1}$,\quad
   $\breve\bA^{({\widetilde m})}  \leftarrow [\breve\bA^{({\widetilde m})},   \bu_{r1}]$,\quad
   $\breve\bA^{({\widetilde m}+1)}  \leftarrow [\breve\bA^{({\widetilde m}+1)},   \bv_{r1}]$\;
\nl   $\bM  \leftarrow \blkdiag(\bM , \1_{1\times J_r})$\;
	}
\;\vspace{-3ex}
{\mtcc{Stage 3b: Rank-$J$ to rank-$R$ Kruskal tensor \dashrule}}
\nl ${\tensor{\widetilde Y}}_R = \llbracket \blambda;  \bA^{(1)}, \ldots, \bA^{({\widetilde m}-1)}, \breve\bA^{({\widetilde m})}, \breve\bA^{({\widetilde m}+1)}, \bA^{({\widetilde m}+1)}, \ldots, \bA^{({\widetilde M})} \rrbracket$\;
  \If{$\sum_{r = 1}^{R}{J_r} >  R$}{
\nl ${\tensor{\widetilde Y}}_J = \llbracket \blambda;  \bA^{(1)} \bM, \ldots, \bA^{({\widetilde m}-1)} \bM, \tilde\bA^{({\widetilde m})}, \tilde\bA^{({\widetilde m}+1)}, \bA^{({\widetilde m}+1)} \bM, \ldots, \bA^{({\widetilde M})} \bM\rrbracket$
\nl ${\llbracket \blambda; \bA^{(1)}, \bA^{(2)}, \ldots, \bA^{(\widetilde M)} \rrbracket}  = \psCPD({\tensor{\widetilde Y}}_J,R, {\tensor{\widetilde Y}}_R)$
}
	\nl $\widetilde M = \widetilde M +1$, 
}
}
}
\;\vspace{-3ex}
 {\mtcc{Stage 4: Refinement  if needed \dashrule}}
\nl $\llbracket \blambda; \bA^{(1)}, \ldots, \bA^{(N)} \rrbracket = \CPD(\tY,R, {\llbracket \blambda; \bA^{(1)}, \bA^{(2)}, \ldots, \bA^{(N)} \rrbracket} )$
}
\BlankLine
\minitab[p{.95\linewidth}]{\hline
$\tucker(\tY, \bsR)$:  rank-$\bsR$ Tucker decomposition of order-$N$ tensor $\tY$ where $\bsR = [R_1, R_2, \ldots, R_N]$.\\
${\tensor{\widehat{Y}}}$ = \CPD($\tY$, $R$, $\tY_{init}$):  approximates an order-$N$ tensor or a tensor in the Kruskal form $\tY$ by a rank-$R$ Kruskal tensor ${\tensor{\widehat{Y}}}$ using initial values $\tY_{init}$.
}
\end{algorithm}
\def\baselinestretch{1.5} 

\subsection{The proposed Algorithm}
 
%

When the tensor $\tY$ is unfolded by a complex unfolding rule $\bl$ which comprises multiple two-modes or $M$-mode unfoldings such as $\bl  = [(1,2),(3,4,5),(7,8)]$, construction of a rank-$J$ structured Kruskal tensor becomes complicated.
In such case, the factor reconstruction process in section~\ref{sec::fold2modes}
or section~\ref{sec::foldMmode} is sequentially applied to mode to be unfolded.
In Algorithm~\ref{alg_3w2Nw_lowrank}, we present a simple version of FCP using low-rank approximation to merge components.
The algorithm reduces the tensor order from $N$ to $M$ (e.g., 3)  specified by the unfolding rule $\bl = [\bl_1, \bl_2, \ldots, \bl_M]$, where each group of modes $\bl_m  = [l_m(1), l_m(2), \ldots, l_m(K_m)]$, $K_m \ge 1$ and $\sum_{m = 1}^{M} K_m = N$. 

Tucker compression can be applied to the unfolded tensor $\tY_{\llbracket \bl \rrbracket}$ \cite{Andersson98}. The factor matrices are sequentially recovered over unfolded modes $m$, i.e., modes have $K_m \ge 2$.
The sequential reconstruction of factors is also executed over $(K_m -1)$ runs. Each run constructs an order-$\widetilde M$ rank-$J$ Kruskal tensor, and approximates it by an order-$\widetilde M$ rank-$R$ Kruskal tensor using $\tensor{\widetilde Y}_{R}$ to initialize. 
It also indicates that ${\tensor{\widehat{Y}}}$ has better approximation error than that of $\tensor{\widetilde Y}_R$.
The tensor order $\widetilde M$ gradually increases from $M$ to $N$.
The full implementation of FCP provided at {\tt{http://www.bsp.brain.riken.jp/{$\sim$}phan/tensorbox.php}} 
includes other multimode unfoldings. 


Although a rank-$R$ CPD of unfolded tensor has lower approximation error than the best rank-$R$ CPD of the original tensor, for difficult data with collinear components or under-rank approximation ($R$ is much lower than the true rank), CPDs of the unfolded tensors and structured Kruskal tensors are often proceeded with a slightly higher rank $R$.

For some cases, a refinement stage may be needed to obtain the best solution. That is the approximation solution after low-rank approximations is used to initialize CPD of the raw data.
This stage often requires lower number of iterations than CPD with random or SVD-based initializations.

 \section{Simulations}\label{sec::simulation}


Throughout the simulations, the ALS algorithm factorized data tensors in 1000 iterations and stopped when $\displaystyle \varepsilon \le 10^{-8}$. The FCP algorithm itself is understood as Algorithm~\ref{alg_3w2Nw_lowrank} with low-rank approximation. Otherwise, the FCP algorithm with rank-one approximation is denoted by R1FCP. ALS was also utilized in FCP to decompose unfolded tensors.

\begin{figure}
\centering
\psfrag{c1c2c3}[t][b]{\scalebox{1}{\color[rgb]{0,0,0}\setlength{\tabcolsep}{0pt}\vspace{3ex}\begin{tabular}{c}\small $(c_1, c_2, c_3)$\end{tabular}}}%
\psfrag{rank-1, direct}[l][l]{\scalebox{1}{\color[rgb]{0,0,0}\setlength{\tabcolsep}{0pt}\begin{tabular}{c}\vspace{.5ex}\tiny rank-1, $[1,2, (3:6)]$\end{tabular}}}%
\psfrag{rank-1, one-mode}[l][l]{\scalebox{1}{\color[rgb]{0,0,0}\setlength{\tabcolsep}{0pt}\begin{tabular}{c}\vspace{.5ex}\tiny rank-1, one-mode\end{tabular}}}%
\psfrag{rank-1, half rule}[l][l]{\scalebox{1}{\color[rgb]{0,0,0}\setlength{\tabcolsep}{0pt}\begin{tabular}{c}\vspace{.5ex}\tiny rank-1, half rule\end{tabular}}}%
\psfrag{low rank, one-mode}[l][l]{\scalebox{1}{\color[rgb]{0,0,0}\setlength{\tabcolsep}{0pt}\begin{tabular}{c}\vspace{.5ex}\tiny low rank, one-mode\end{tabular}}}%
\psfrag{low rank, half rule}[l][l]{\scalebox{1}{\color[rgb]{0,0,0}\setlength{\tabcolsep}{0pt}\begin{tabular}{c}\vspace{.5ex}\tiny low rank, half-rule\end{tabular}}}%
\psfrag{rank-1, ``good folding''}[l][l]{\scalebox{.8}{\color[rgb]{0,0,0}\setlength{\tabcolsep}{0pt}\vspace{3ex}\begin{tabular}{c}\footnotesize Alg.~\ref{alg_3w2Nw_rank1}, good folding\end{tabular}}}%
\psfrag{rank-1, ``bad folding''}[l][l]{\scalebox{.8}{\color[rgb]{0,0,0}\setlength{\tabcolsep}{0pt}\vspace{3ex}\begin{tabular}{c}\footnotesize Alg.~\ref{alg_3w2Nw_rank1}, ``bad folding''\end{tabular}}}%
\psfrag{low rank, ``good folding''}[l][l]{\scalebox{.8}{\color[rgb]{0,0,0}\setlength{\tabcolsep}{0pt}\vspace{3ex}\begin{tabular}{c}\footnotesize Alg.~\ref{alg_3w2Nw_lowrank}, ``good folding''\end{tabular}}}%
\psfrag{low rank, ``bad folding''}[l][l]{\scalebox{.8}{\color[rgb]{0,0,0}\setlength{\tabcolsep}{0pt}\vspace{3ex}\begin{tabular}{c}\footnotesize Alg.~\ref{alg_3w2Nw_lowrank}, ``bad folding''\end{tabular}}}%
\psfrag{ALS}[l][l]{\scalebox{.8}{\color[rgb]{0,0,0}\setlength{\tabcolsep}{0pt}\vspace{3ex}\begin{tabular}{c}\footnotesize ALS\end{tabular}}}%
\subfigure[SAE loss.]
{
\includegraphics[width=.48\linewidth, trim = 0.0cm 0cm 0cm 1cm,clip=true]{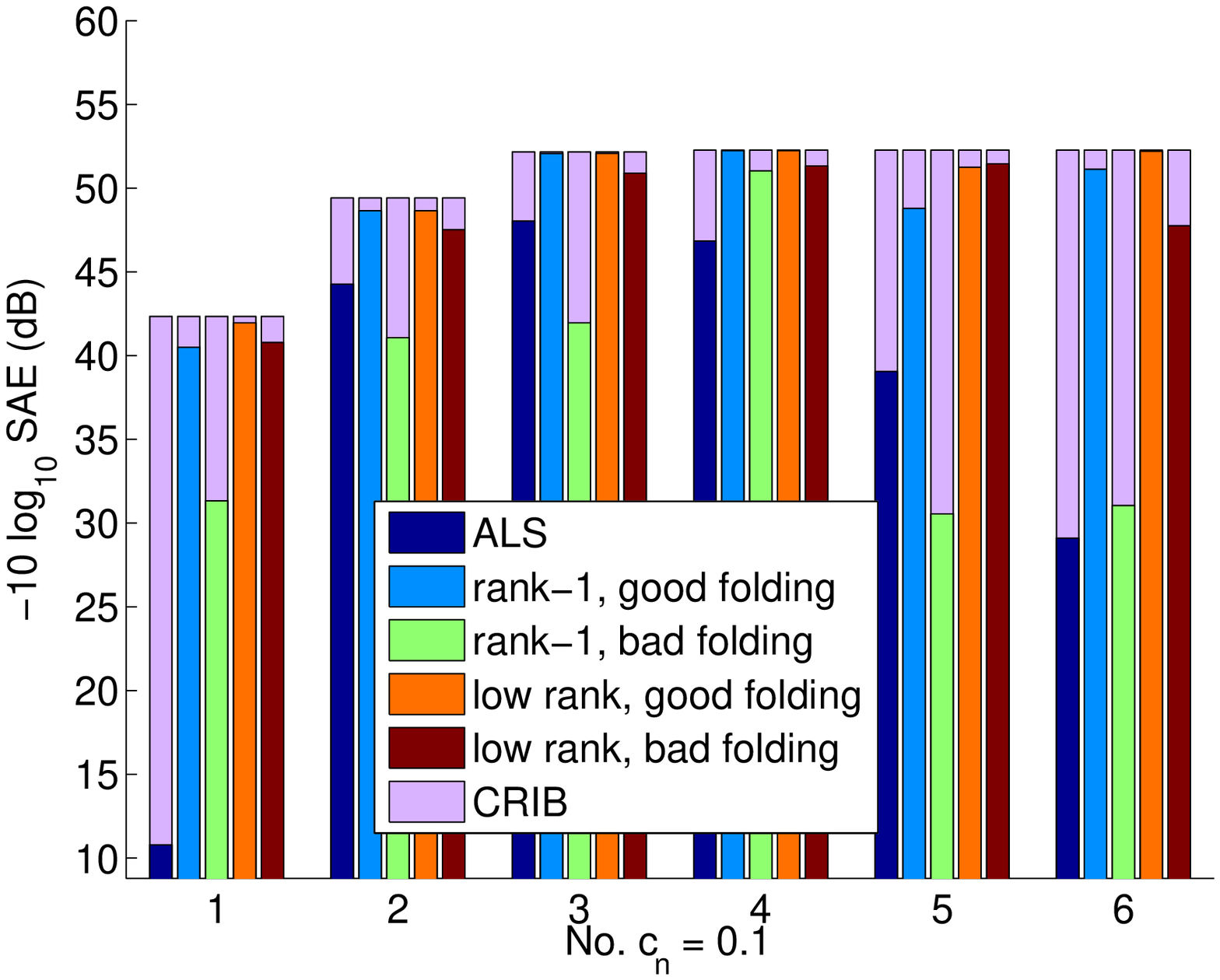}\label{fig_mLOSS_N6I20R20_SNR0_newrule}
}\hfill
\subfigure[Execution time.]
{
\includegraphics[width=.48\linewidth, trim = 0.0cm 0cm 0cm 1cm,clip=true]{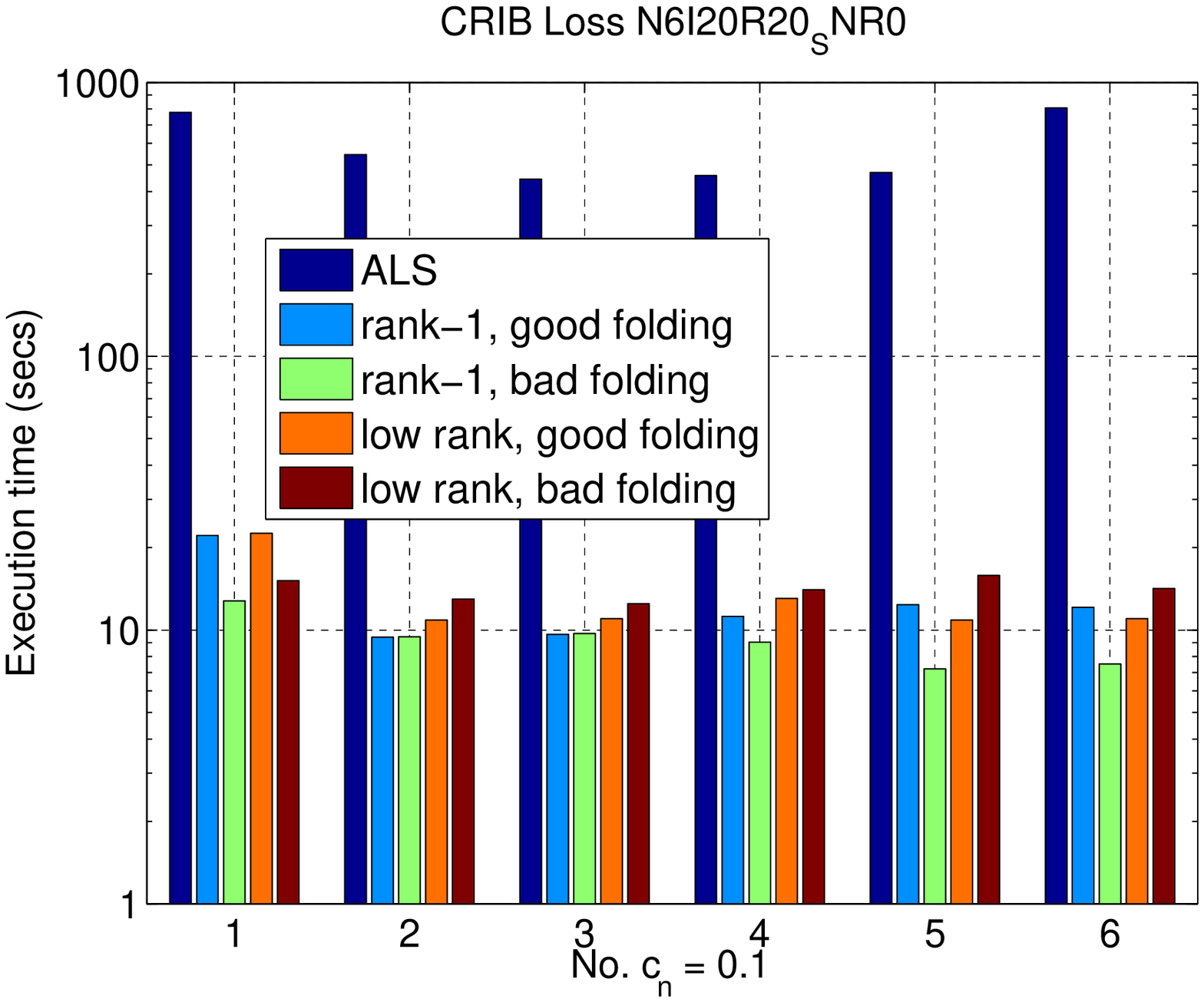}\label{fig_Exectime_N6I20R20_SNR0_newrule}}
\caption{Illustration of MSAE loss and execution time  averaged over 30 MC runs in decomposition of order-6 rank-20 tensors of size $I = R = 20$ corrupted with additive Gaussian noise of 0 dB SNR, and $c_n \in \{0.1, 0.95\}$.}
\label{fig_ex_N6R20_c1c2c3}
\end{figure}

\begin{example}{\textbf{(Decomposition of order-6 tensors.)}}\label{ex_7}
\end{example}
We analyzed the mean SAEs (MSAE) of algorithms in decomposition of order-6 tensors with size $I_n = R = 20$ by varying the number of low collinear factors from 1 to 6 with $[c_1 \le c_2 \le \cdots \le c_6]$. 
Tensors were corrupted with additive Gaussian noise of 0 dB SNR. 
ALS was not efficient and often got stuck in local minima.
MSAEs over all the estimated components by ALS were clearly lower than CRIB, especially when there were 5 collinear factors (the first test case in Fig~\ref{fig_mLOSS_N6I20R20_SNR0_newrule}).
The FCP method was executed with ``good unfolding'' rules suggested by the strategy in Section~\ref{sec::unfoldingstrategy} and ``bad'' ones which violated the unfolding strategy listed in Table~\ref{tab_6D}.

\begin{table}
\centering
\footnotesize
\caption{Comparison of MSAEs (in dB) of ALS and FCP with CRIB in decomposition of order-6 rank-20 tensors of size $I_n = 20$, for all $n$. Correlation coefficients of factors $c_n$ are  0.1 or 0.9 and $[c_1 \le  c_2 \le \cdots \le c_6]$. The performance was measured by varying the number of $c_n = 0.1$.
R1FCP ($J_r = 1$, for all $r$) was sensitive to unfolding rules.}
\label{tab_6D}
\begin{tabular}{llllllllll}
& & & \multicolumn{3}{c}{``Good unfolding''} & \multicolumn{3}{c}{``Bad unfolding''}\\\cline{4-6}\cline{7-9}
No. $0.1$		& CRIB &	ALS &   Unfolding  & 	Alg.~\ref{alg_3w2Nw_rank1}  & Alg.~\ref{alg_3w2Nw_lowrank} & Unfolding   & 	Alg.~\ref{alg_3w2Nw_rank1}  & Alg.~\ref{alg_3w2Nw_lowrank}\\\hline
1 &	42.33&10.8  &$[1, (2, 3, 4), (5, 6)]$ 	&		40.49	  &41.95		  &$[2, 3, (1, 4, 5, 6)]$		&{\color{red}{\bf{31.33}}}	 &40.78\\
2 & 	49.40& 44.3 &$[1, 2, (3, 4, 5, 6)]$    	&		48.64	  &48.65		  &$[3, 4, (1, 2, 5, 6)]$		&{\color{red}{\bf{41.06}}}	 &47.50\\
3 &	52.16& 48.0 &$[1, 2, (3, 4, 5, 6)]$	&		52.06	  &52.07		  &$[4, 5, (1, 2, 3, 6)]$		&{\color{red}{\bf{41.94}}} &50.89\\
4 &	52.26& 46.8 &$[1, 2, (3, 4, 5, 6)]$	&		52.23	  &52.23		  &$[1, (2, 3),(4, 5, 6)]$	&51.02	 &51.31\\
5 &	52.26& 36.2 &$[1, 2, (3, 4, 5, 6)]$	&		48.78	  &51.24		  &$[(1, 2), (3, 4),(5, 6)]$	&{\color{red}{\bf{30.54}}} &51.44\\
6 &	52.26& 30.5 &$[1, 2, (3, 4, 5, 6)]$	&		51.13	  &52.19		  &$[(1, 2), (3, 4),(5, 6)]$	&{\color{red}{\bf{31.04}}}	 &47.75
\end{tabular}
\end{table}
Performance of R1FCP (Algorithm~\ref{alg_3w2Nw_rank1}) was strongly affected by the unfolding rules. Its SAE loss was up to 21dB with ``bad unfoldings''.
For all the test cases, FCP with low-rank approximations (i.e., Algorithm~\ref{alg_3w2Nw_lowrank}) obtained high performance even with ``bad unfolding'' rules.
In addition, FCP was much faster than ALS. 
FCP factorized order-6 tensors in 10-20 seconds while ALS completed the similar tasks in 500-800 seconds. 
Finally, in this simulation, FCP was 47 times faster on average than ALS. 

\begin{example}{\textbf{(Factorization of Event-Related EEG Time-Frequency Representation.)}}\label{ex_8}
\end{example} 
This example illustrates an application of CPD for analysis of  real-world EEG data \cite{Morup,Morup06} which consisted of 28 inter-trial phase coherence (ITPC) measurements \cite{Tallon1996} of EEG signals of 14 subjects during a proprioceptive pull of the left and right hands.
The whole ITPC data set was organized as a 4-way tensor of
28 measurements $\times$ 61 frequency bins $\times$  64 channels $\times$ 72 time frames.
The first 14 measurements were associated to a group of the left hand stimuli, while the other ones were with the right hand stimuli.
The order-4 ITPC tensor can be fitted by a multilinear CP model. 
M{\o}rup \emph{et al.} analyzed the dataset by nonnegative CP of three components and Tucker components and compared them with components extracted by NMF and ICA \cite{Morup06}.


%
%
%
%
%

\begin{table}
\caption{Comparison of fit (\%) values in factorization of the ITPC tensor by ALS and FCP in Example~\ref{ex_8}. Rank-1 FCP (i.e., Algorithm~\ref{alg_3w2Nw_rank1}) completely failed in this example. Strikethrough values mean that the algorithm did not converge to the desired solution.}
\label{tab_fit_EEGpull}
\centering
\begin{tabular*}{1\textwidth}{@{\extracolsep{\fill}}c@{\hspace{1ex}}c@{\hspace{1ex}}c@{\hspace{1ex}}c@{\hspace{1ex}}c@{\hspace{.5ex}}c@{\hspace{.5ex}}c@{}}
R 	&ALS           &  \multicolumn{2}{c}{Rank-1 FCP} & Low-rank FCP   &Tucker $\rightarrow$ ALS        & Tucker $\rightarrow$ Alg.~\ref{alg_3w2Nw_lowrank}\\
 & & [(1,2), 3, 4] & [1,2, (3, 4)]
\\\hline
5	&37.7 $\pm$ 0.17  & 32.7   &36.0  &37.7  $\pm$ 0.17   &36.2   $\pm$0.03    &36.9   $\pm$0.00  \\
8	&43.8 $\pm$ 0.13  & 26.7    &42.1	&43.8  $\pm$ 0.00   &42.7   $\pm$0.01    &42.9   $\pm$0.58  \\
11	&48.0 $\pm$ 0.04  & 19.5    &32.5	&47.9  $\pm$ 0.16   &47.1   $\pm$0.09    &47.5   $\pm$0.08  \\
20	&56.0 $\pm$ 0.12   & \st{-61.1}  &\st{-10.8}  	&55.9  $\pm$ 0.16   &55.5   $\pm$0.18    &55.8   $\pm$0.13  \\
30	&61.1 $\pm$ 0.09   & \st{-519}   &\st{-301.0}	&61.1  $\pm$ 0.09   &60.8   $\pm$0.07    &60.9   $\pm$0.13  \\
40	&64.5 $\pm$ 0.08   & \st{-649}   &\st{-319.1}	&64.4  $\pm$ 0.14   &64.2   $\pm$0.08    &64.3   $\pm$0.18  \\
60	&68.7 $\pm$ 0.02   & \st{-1295}& \st{-432.8} 	&68.1  $\pm$ 0.05   &68.6   $\pm$0.09    &68.6   $\pm$0.10 \\ 
72	&70.4 $\pm$ 0.02   & \st{-7384}& \st{-535.0}   	&69.9  $\pm$ 0.08   \\ \hline
\end{tabular*}
\end{table}

\begin{table}[t]
\centering
\caption{Performance of rank-1 FCP with different unfolding rules in decomposition of the ITPC tensor in Example~\ref{ex_8}. Strikethrough values mean that the algorithm did not converge to the desired solution.}\label{tab_r1FCP_ITPC}
\begin{tabular}{@{\extracolsep{\fill}}cc@{\hspace{2ex}}c@{\hspace{1ex}}c@{\hspace{2ex}}ccc@{\hspace{1ex}}c@{\hspace{1ex}}c@{}}
\multirow{3}{*}{$R$} & \multicolumn{4}{c}{Rank-1 FCP} & &\multicolumn{2}{c}{ALS} \\\cline{2-5}\cline{7-8}
 & \multirow{2}{*}{\minitab[c]{Unfolding\\ rules}}	&Fit 		&Time 	&Collinearity degree  && Fit & Time\\
&&(\%)& (secs) &$[c_1, c_2, c_3, c_4]$ &&  (\%) & (secs) \\\cline{1-5}\cline{7-8}
\multirow{4}{*}{8} 
 &[(1, 2), 3, 4]  	&26.7 	&2.15     	&[0.48, 0.70, 0.54, 0.89]	&& \multirow{4}{*}{43.8} & \multirow{4}{*}{242} & \multirow{4}{*}{} \\
 &[1, (2, 3), 4]	&29.9	&5.72	&[0.37, 0.50, 0.40, 0.83]	&\\
 &[1, 2, (3, 4)]	&42.1	&5.59	&[0.34, 0.50, 0.49, 0.72]	 &\\
 &[1, 3, (2, 4)]	&41.3	&4.67	&[0.49, 0.52, 0.50, 0.75] 	&\\\cline{1-5}\cline{7-8}
\multirow{4}{*}{20} 
&[(1, 2), 3, 4]	&\st{-61.1}   	&5.79     	&[0.43, 0.62, 0.50, 0.89] && \multirow{4}{*}{56.0} & \multirow{4}{*}{752} &  \\
&[1, (2, 3), 4]	&28.3	&12.12	&[0.34, 0.45, 0.56, 0.90]	\\
&[1, 2, (3, 4)] 	&15.4	&10.39	&[0.33, 0.48, 0.46, 0.86] \\
&[1, 3, (2, 4)]	&\st{-22.6} 	&9.17	&[0.39, 0.57, 0.69, 0.92] \\
\cline{1-5}\cline{7-8}
\end{tabular}
\end{table}

  \begin{figure}
\centering
\psfrag{3W2NW}[cl][cl]{\scalebox{.9}{\color[rgb]{0,0,0}\setlength{\tabcolsep}{0pt}\begin{tabular}{c}\footnotesize FCP\end{tabular}}}%
\psfrag{Tucker->3W2NW}[cl][cl]{\scalebox{.9}{\color[rgb]{0,0,0}\setlength{\tabcolsep}{0pt}\begin{tabular}{c}\footnotesize Tucker$\rightarrow$FCP\end{tabular}}}%
\psfrag{Tucker->ALS}[cl][cl]{\scalebox{.9}{\color[rgb]{0,0,0}\setlength{\tabcolsep}{0pt}\begin{tabular}{c}\footnotesize Tucker$\rightarrow$ALS\end{tabular}}}%
\psfrag{ALS}[cl][cl]{\scalebox{.9}{\color[rgb]{0,0,0}\setlength{\tabcolsep}{0pt}\begin{tabular}{c}\footnotesize ALS\end{tabular}}}%
\subfigure[Execution time (seconds)]
{
\includegraphics[width=.48\linewidth, trim = 0.0cm 0cm 0cm 0cm,clip=true]{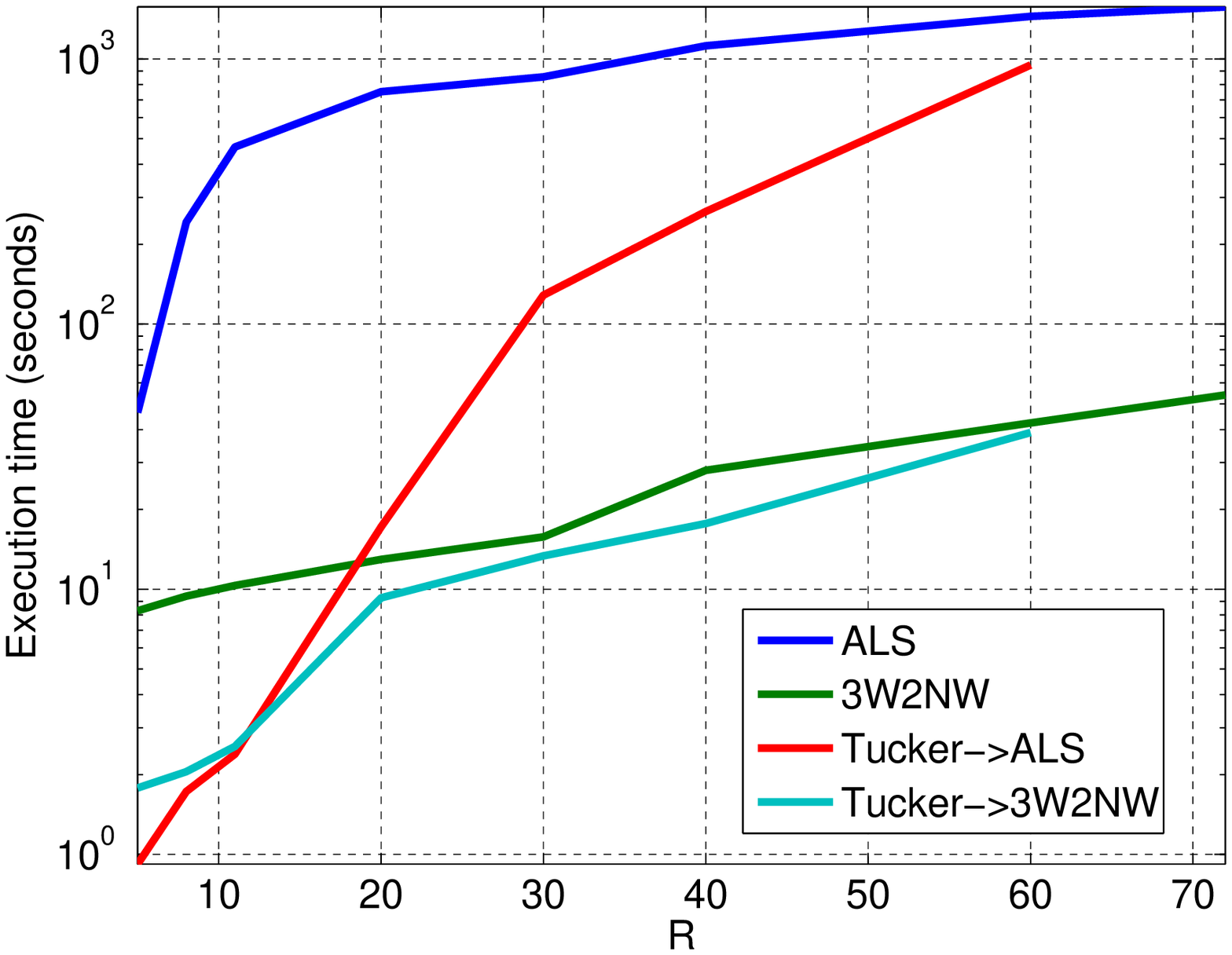}\label{fig_EEGpull_Exectime}}
\hfill
\subfigure[Execution time as $R = 20$]
{
\includegraphics[width=.48\linewidth, trim = 0.0cm 0cm 0cm 0cm,clip=true]{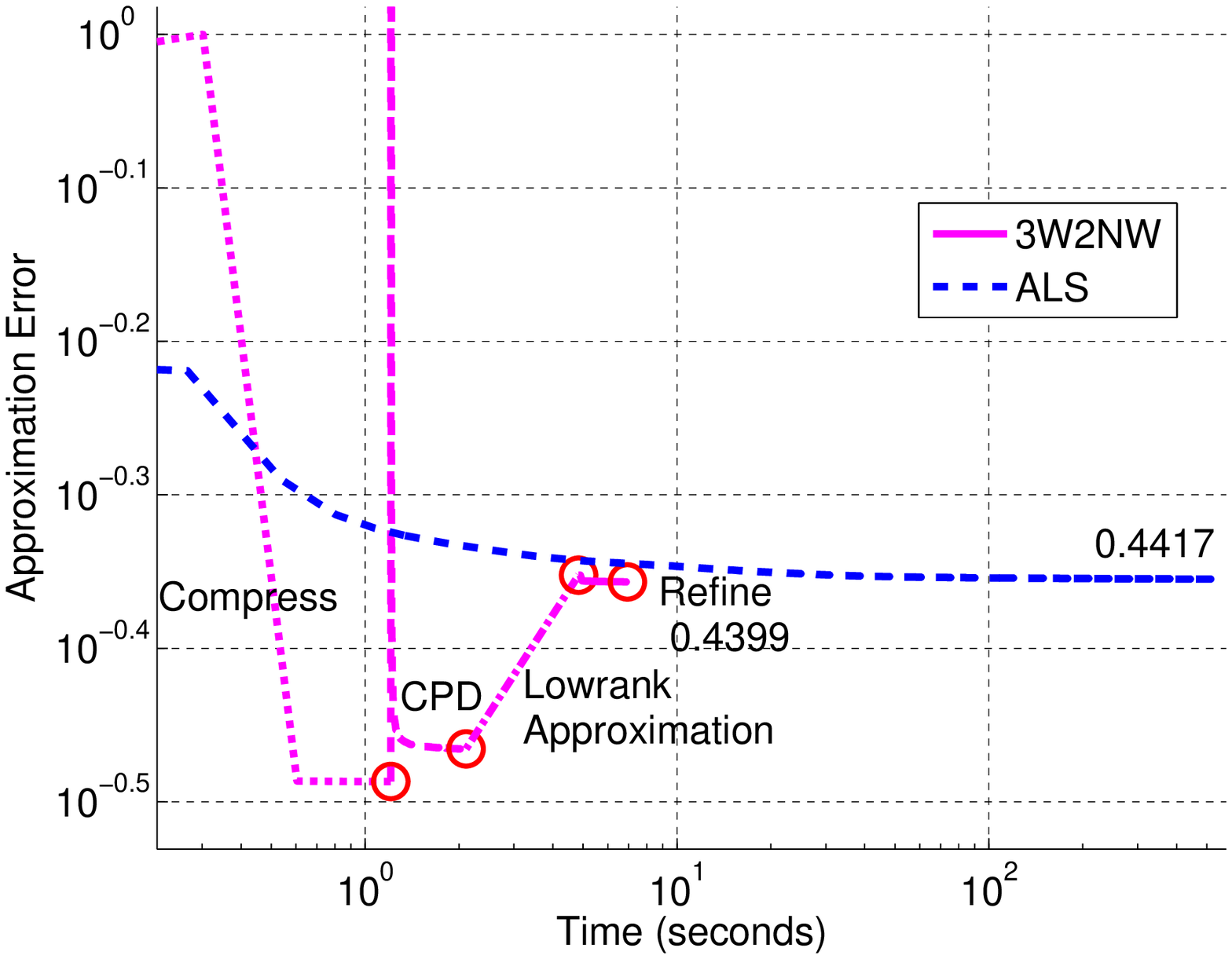}\label{fig_EEGpull_errvstime}}
\caption{Illustration of execution times (seconds) of ALS and FCP for factorization of order-4 ITPC tensor in Example~\ref{ex_8}.}
\label{fig_ex_EEG}
\end{figure}

   \begin{figure}
\centering
\subfigure[$R = 8$]{
\includegraphics[width=.47\linewidth, trim = 0.0cm 0cm 0cm 0cm,clip=true]{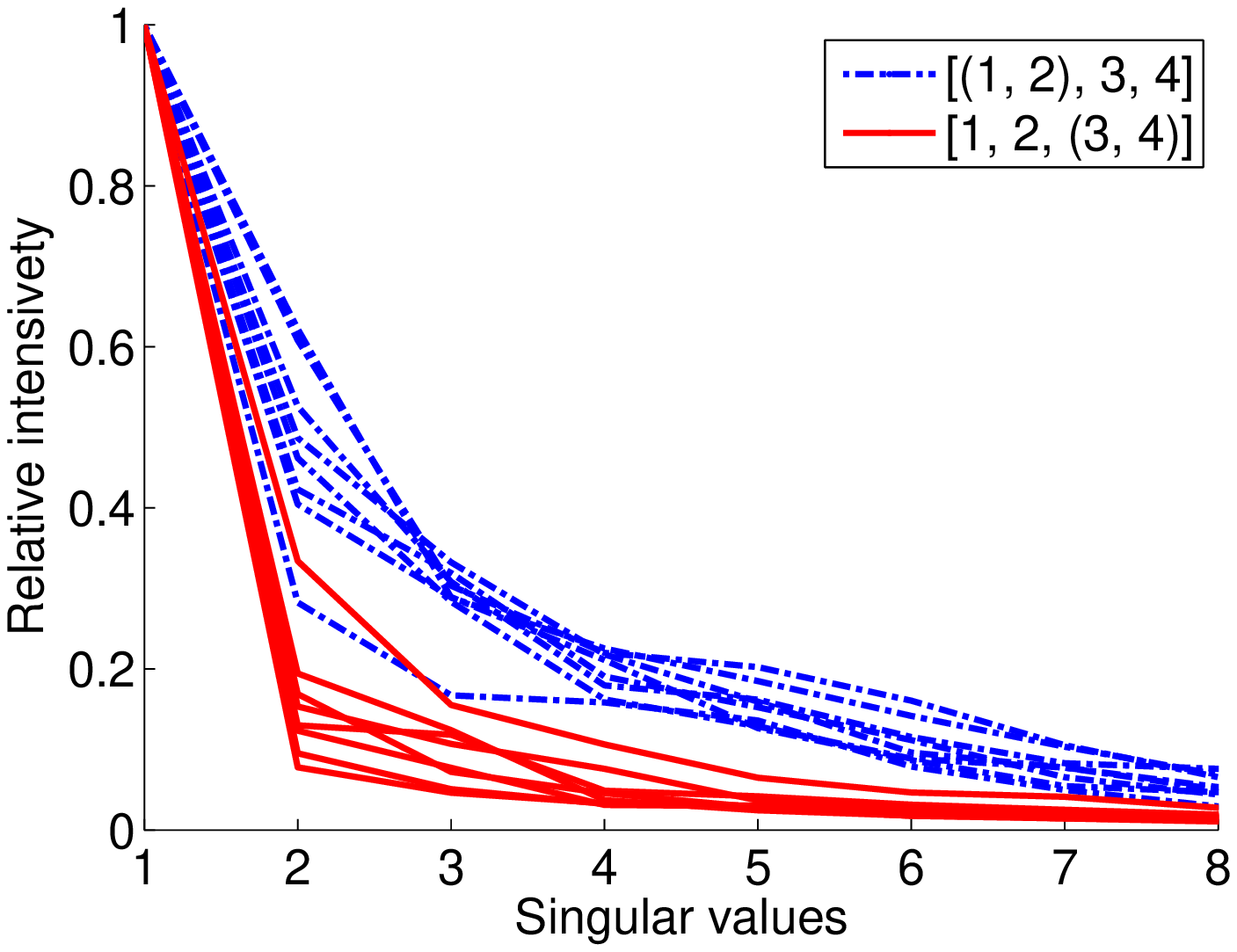}\label{fig_eegpull_svd_r1fcp_R8_28616472__12_vs_34_}}
\hfill
\subfigure[$R = 20$]{
\includegraphics[width=.47\linewidth, trim = 0.0cm 0cm 0cm 0cm,clip=true]{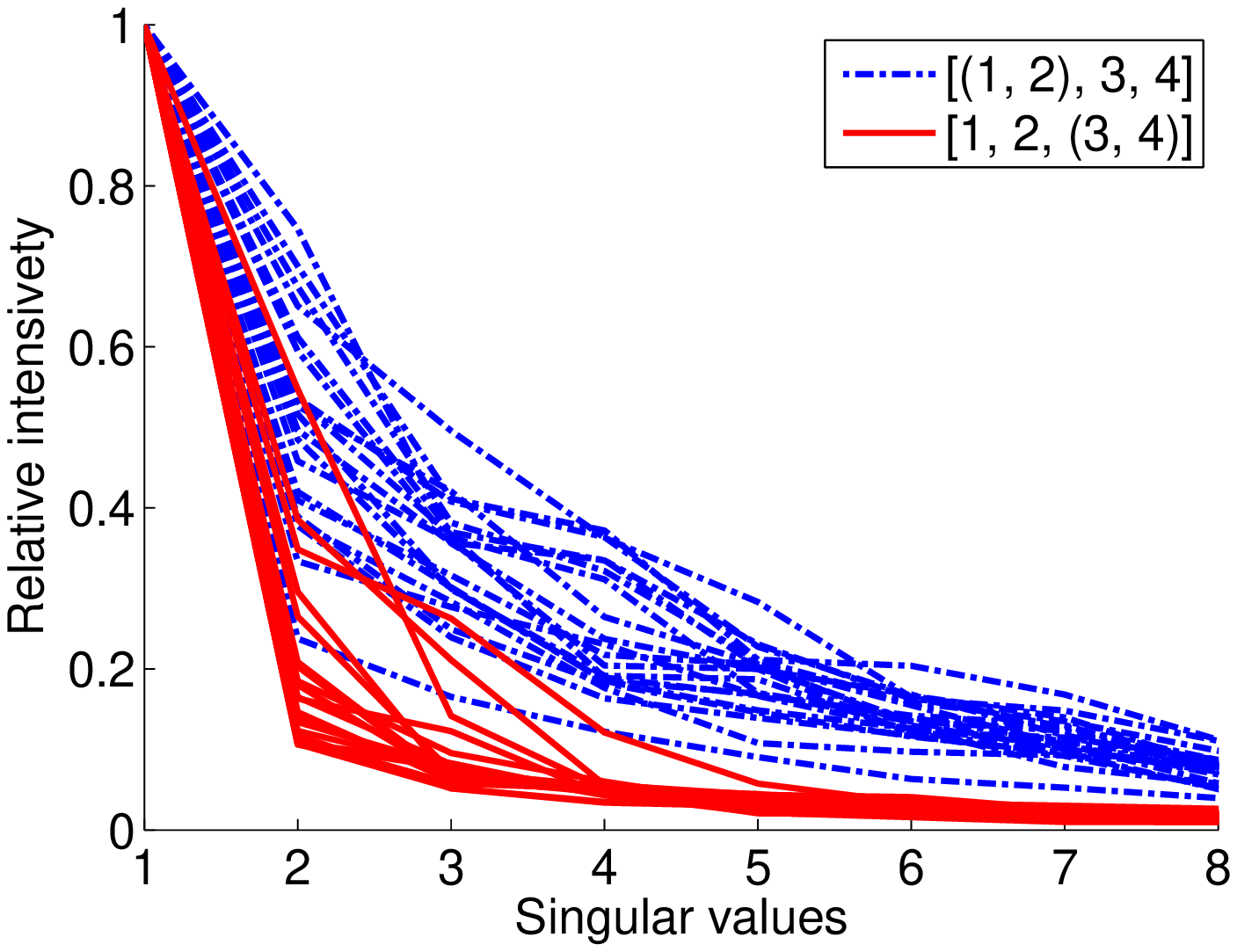}\label{fig_eegpull_svd_r1fcp_R20_28616472__12_vs_34_}}
\caption{Illustration of leading singular values of matrices $\bF_r$ ($r = 1, 2, \ldots, R$) reshaped from components estimated from the ITPC tensor with different unfolding rules $[(1,2),3,4]$ and $[1,2,(3,4)]$. Rank-1 FCP (i.e., Algorithm~\ref{alg_3w2Nw_rank1}) failed in this experiment because this algorithm worked only if all $\bF_r$ were rank-one matrices.}
\label{fig_ex_EEG_svd}
\end{figure}

In this example, our aim was to compare the factorization time of ALS and FCP over various $R$  in the range of [5, 72] with and without a Tucker compression prior to the CP decompositions. 
The FCP method employed ALS to factorize the order-3 unfolded tensor, and the fast ALS for the structured Kruskal tensors (see in Appendix \ref{sec::app_fastCP}).
Interpretation of the results can be found in \cite{Morup,Morup06}.
The low-rank FCP algorithm was applied with the unfolding rule $\bl = [1, 2, (3, 4)]$.

Execution time for each algorithm was averaged over 10 Monte Carlo runs with different initial values and illustrated in Fig.~\ref{fig_EEGpull_Exectime} for various $R$.
For relatively low rank $R$, a prior Tucker compression sped up ALS, and made it more efficient than FCP when $R \le 11$. 
The reason is explained by compression time for unfolding tensor in FCP.
However, this acceleration technique was less efficient as $R \rightarrow I_n$ and inapplicable to ALS for $R\ge I_n$.
FCP significantly reduced the execution time of ALS by a factor of 5-60 times, and was slightly improved by the prior compression.
Comparison of fits explained by algorithms in Table~\ref{tab_fit_EEGpull} indicates that while FCP with Algorithm~\ref{alg_3w2Nw_lowrank} quickly factorized the data, it still maintained fit equivalent to ALS.

For this data, the rank-one FCP algorithm (i.e., Algorithm~\ref{alg_3w2Nw_rank1}), unfortunately, did not work well. Fits of this algorithm are given in Table~\ref{tab_fit_EEGpull}.
Performance of this algorithm with several unfolding rules including $[(1,2), 3, 4]$, $[1, (2, 3), 4]$, $[1,2, (3, 4)]$ and $[1,3, (2, 4)]$ is compared in Table~\ref{tab_r1FCP_ITPC}.
When $R = 8$ and using the rule $\bl = [(1,2), 3, 4]$, the rank-1 FCP algorithm showed the worst performance with a fit of $26.7\%$ which was not competitive to a fit of 43.8\% obtained by ALS.
The average collinearity degrees of the estimated components  $c_n = [0.48, {\bf0.70}, 0.54, {\bf0.89}]$ indicates that we should not fold modes 1 and 2; in addition, folding modes 2 and 4 which had the largest collinear degrees is suggested, i.e., the unfolding rule $\bl = [1, 3, (2, 4)]$.
It is clear to see that the unfolding rule $\bl = [1, 3, (2, 4)]$ significantly improved performance of the rank-1 FCP algorithm with a fit of 41.3\%. Moreover, the unfolding rule $\bl = [1, 3, (2, 4)]$  was also suggested according to the average collinear degrees $c_n = [0.37, {\bf 0.50}, 0.40, {\bf 0.83}]$ achieved when applying the unfolding rule $\bl = [1, (2,3), 4]$.
This confirms the applicability of the suggested unfolding strategy.
For this test case, the unfolding rule $\bl = [1, 2, (3, 4)]$ allowed to achieve the best fit of $42.1\%$, although this rule was not suggested by the strategy. This can be understood due to the fact that the average collinear degrees of modes 2 and 3 were very similar (0.50 and 0.49, or 0.52 and 0.50, see in Table~\ref{tab_r1FCP_ITPC}).

For higher ranks, e.g.,  $R \ge 11$, FCP with rank-one approximation completely failed. The unfolding strategy did not help anymore (see fit values in Table~\ref{tab_fit_EEGpull}).
In Fig.~\ref{fig_ex_EEG_svd}, we display leading singular values of reshaped matrices $\bF_{r}$ ($r = 1, 2, \ldots, R$) from the estimated components for $R = 8$ and 20.
The results indicate that $\bF_r$ were not rank-one matrices, especially the matrices received when using the rule $\bl = [(1, 2), 3, 4]$. Note that the rank-one FCP algorithm works if and only if all $\bF_r$ are rank-one.
This also confirms that the low-rank FCP algorithm was appropriate for this data.


Fig.~\ref{fig_EEGpull_errvstime} illustrates the relative approximation errors $\varepsilon = \frac{\|\tY - \tensor{\hat Y}\|_F}{\|\tY\|_F}$ of ALS and FCP for $R = 20$ as functions of the execution time.
ALS took 536.5 seconds to converge. FCP took 1.2 seconds for compression, 0.9 seconds for CPD of the order-3 unfolded tensor, 2.73 seconds for low-rank approximations, 2.1 seconds for the refinement stage. ALS and FCP converged to the relative approximation errors $\varepsilon_{ALS}$ = 0.4417, while $\varepsilon_{FCP}$ = 0.4399, respectively.

\begin{example}{\textbf{(Decomposition of Gabor tensor of the ORL face database.)}}\label{ex_10}
\end{example} 
%
This example illustrated classification of the ORL face database \cite{Samaria1994} 
consisting of 400 faces for 40 subjects.
We constructed Gabor feature tensors for 8 different
orientations at 4 scales which were then down-sampled to 16 $\times$ 16 $\times$ 8  $\times$ 4 $\times$ 400 dimensional tensor $\tY$.
The unfolding $\bl = [1, 2, (3, 4,5)]$ was applied to unfold $\tY$ to be an order-3 tensor. 
ALS \cite{Nwaytoolbox} factorized both $\tY$ and $\tY_{\llbracket \bl \rrbracket}$ into $R = 30, 40, 60$ rank-1 tensors
in 1000 iterations, and stopped when $\displaystyle {\|\varepsilon - \varepsilon_{old}\|} \le 10^{-6} \varepsilon$ where $\varepsilon = \|\tY - \hat{\tY}\|_F^2$. 
The rank-one FCP algorithm did not work for this data.
For example, when $R = 10$ and applying the rule $\bl = [1, 2, (3, 4,5)]$, R1FCP explained the data with a fit of -31.2\%, and yielded average collinearity degrees of $[0.60, 0.66, 0.64, 0.95, 0.64]$. Although a further decomposition with the unfolding rule $\bl = [1, (2,3), (4,5)]$ achieved a better fit of 44.8\%, this result was much lower than a fit of $ 54.5\%$ obtained by ALS and FCP. 



The factor $\bA^{(5)} \in \Real^{400 \times R}$ comprised compressed features which were used to cluster faces using the K-means algorithm.
Table~\ref{tab_orl} compares performance of two algorithms including execution time, fit, accuracy (ACC \%) and normalized mutual information (NMI).
For $R = 40$, ALS factorized $\tY$ in 1599 seconds while FCP completed the task in only 39 seconds 
with a slightly reduction of fit ($\approx 0.17\%$). 
For $R = 60$, ALS was extremely time consuming, required 16962 seconds while FCP only took 105 seconds.
Regarding the clustering accuracy, features extracted by FCP still achieved comparable performance as those obtained by ALS.

\begin{table}
\centering
\small
\caption{Comparison between ALS and FCP (Alg.~\ref{alg_3w2Nw_lowrank}) in factorization of order-5 Gabor tensor constructed from the ORL face dataset.}\label{tab_orl}
\begin{tabular}{@{}l@{\hspace{1ex}}cc@{\hspace{1ex}}r@{\hspace{1ex}}c@{}cc@{}}
$R$   & \minitab[c]{\parbox{3em}{Algorithm}} &  Fit (\%) & \minitab[c]{Time\\(seconds)}  & \minitab[c]{Ratio $\frac{\mbox{ALS}}{\mbox{FCP}}$}  & ACC ($\%$)  & NMI ($\%$)\\
\hline
\multirow{2}{*}{30} & FCP    &60.59  	 &24 	  &\multirow{2}{*}{\bf39}  & 85.00 	&  92.91   	\\
			   & ALS      &60.56 	 &927  &  & 85.25 	&  93.22		\\
\hline
\multirow{2}{*}{40} & FCP   & 62.46 	 &39 	    &\multirow{2}{*}{\bf41}  & 84.25 	 & 92.57   	\\
			   & ALS     & 62.63  	 &1599 &     & 85.50  & 93.68		\\
\hline
\multirow{2}{*}{60} & FCP    & 65.47      &105   &\multirow{2}{*}{\bf162}  & 83.00 	  &91.62   	\\
			   & ALS      & 65.64 	 &16962 &   & 81.38 	  &91.44		\\
\hline	   
\end{tabular}
\end{table}
  



\section{Conclusions}\label{sec:conclusion}

The fast algorithm has been proposed for high order and relatively large scale CPD. The method decomposes the unfolded tensor in lower dimension which is often of size $R \times R \times R$ instead of the original data tensor. Higher order structured Kruskal tensors are then generated, and approximated by rank-$R$ tensors in the Kruskal form using the fast algorithms for structured CPD. Efficiency of the strategy proposed for tensor unfoldings has been proven for real-world data.
In addition, one important conclusion drawn from our study is that factor matrices in orthogonally constrained CPD for high order tensor can be estimated through decomposition of order-3 unfolded tensor without any loss of accuracy.
Finally, the proposed FCP algorithm has been shown 40-160 times faster than ALS for decomposition of order-5 and -6 tensors.

\appendices

\section{Algorithms for structured CPD}\label{sec::app_fastCP}
\subsection{Low complexity CP gradients}
CP gradients of $\tensor{\widetilde Y}$ in Lemma~(\ref{lem_rankJ_CP}) with respect to $\bA^{(n)}$ is quickly computed without construction of $\tensor{\widetilde Y}$ as follows
\be
	\left(\tilde\bY_{(n)} - \hat\bY_{(n)}\right) \, \left(\bigodot_{k\neq n} \bA^{(k)} \right)  = {\tilde\bA^{(n)}} \, \diag(\tilde\blambda) \, \bW_n -  {\bA^{(n)}} \, \bGamma_n  \label{equ_cpgradient_rankJ}
\ee
where $\bW_n = \left( \bigcircledast_{k\neq n} (\tilde\bA^{(n) T}  \, \bA^{(n)}) \right)$ and $\bGamma_n = \left( \bigcircledast_{k\neq n} (\bA^{(k) T}  \, \bA^{(k)}) \right)$.
%
By replacing $\tilde\bA^{(n)}  = \bB^{(n)} \bM$, and employing  
$\left(\bM^T \bA\right) \circledast \left(\bM^T\bB\right) = \bM\,^T \left(\bA\circledast \bB\right)$, $\bM\, \left(\left(\bM^T \bA\right) \circledast \bB\right) = \bA \circledast \left(\bM \bB\right)$,
the first term in (\ref{equ_cpgradient_rankJ}) is further expressed for $n = 1, \ldots, N-2$ by
\begin{align}
{\tilde\bA^{(n)}} \, \diag(\tilde\blambda) \, \bW_{n}
	&= \bB^{(n)}  \, \left(\left(\bigcircledast_{k = 1\hfill \atop k\neq n\hfill}^{N-2} (\bB^{(k) T}  \, \bA^{(k)})\right)  \circledast  \bK
	\right), \notag \\
 \bK &= 	\left[\begin{array}{@{}c@{}}
	\lambda_1	\1_{R}^T \left(	\left(\bU_{1}^T \, \bA^{(N-1)} \right) \circledast \left(\widetilde\bV_{r}^T \, \bA^{(N)} \right)\right)\\
		\vdots\\
	\lambda_R \1_{R}^T \left(	\left(\bU_{R}^T \, \bA^{(N-1)} \right) \circledast \left(\widetilde\bV_{R}^T \, \bA^{(N)} \right)\right)\\
	\end{array}\right] \notag
\end{align}
 and
\be
{\tilde\bA^{(N-1)}} \, \diag(\tilde\blambda) \, \bW_{N-1}
 &=& {\tilde\bA^{(N-1)}}  \, \left[
 \begin{array}{c}
\lambda_1 \widetilde\bV^{T}_1 \, \bA^{(N)}  \, \diag(\bomega_1) \\
\vdots\\
\lambda_R \widetilde\bV^{T}_R \, \bA^{(N)}  \, \diag(\bomega_R)
 \end{array}
 \right] , \notag \\
{\tilde\bA^{(N)}} \, \diag(\tilde\blambda) \, \bW_N
 &=& {\tilde\bA^{(N)}}  \, \left[
 \begin{array}{c}
\lambda_1 \bU^{T}_1 \, \bA^{(N-1)}  \, \diag(\bomega_1) \\
\vdots\\
\lambda_R \bU^{T}_R \, \bA^{(N-1)}  \, \diag(\bomega_R)
 \end{array}
 \right]\, ,\notag
 	\ee	
where $\bOmega = [\bomega_{r}] = \bigcircledast_{\scriptsize k = 1}^{N-2} (\bA^{(k) T}  \, \bB^{(k)})$. 	
For each $n = 1, 2, \ldots, N$, such computation has a low computational complexity of order $\displaystyle \mathcal O\left(R^2 \left(N + \sum_{n = 1}^{N-2}{I_n} \right) + JR \left(I_{N-1} + I_N\right)\right)$.
For tensors which have $I_n = I$, for all $n$, in the worst case, $J_r \approx  I$, then $J \approx R I$, we have
\be
\displaystyle \mathcal O\left(R^2 \left(N + \sum_{n = 1}^{N-2}{I_n} \right) + JR \left(I_{N-1} + I_N\right)\right) 
= \mathcal O\left(NR^2 I^2 \right)  
\ll  \mathcal O\left(R I^{N}\right), \notag
\ee
for $ N \ge4$.

\subsection{Fast algorithms for structured CPD}
By employing the fast CP gradient in previous section, most CP algorithms can be rewritten to estimate $\bA^{(n)}$ from the structured tensors $\tensor{\widetilde Y}$ in Lemma~(\ref{lem_rankJ_CP}). 
For example, the ALS algorithm is given by
\be
	\bA^{(n)}& \leftarrow& \bB^{(n)}  \, \left(\left(\bigcircledast_{k = 1\hfill \atop k\neq n\hfill}^{N-2} (\bB^{(n) T} \bA^{(n)})\right)  \circledast 
\bK
	\right)  \,  \bGamma_n^{-1}\,, \quad n = 1, \ldots, N-2.\notag
\ee

\section{Cram\'er-Rao induced bound for angular error}\label{sec::CRIB}

Denote by $\alpha_1$ the mutual angle between the true factor $\ba^{(1)}_1$ and its estimate ${\hat\ba}^{(1)}_1$
\begin{equation}
\alpha_1 = \mbox{acos}\frac{\ba^{(1)T}_1{\hat\ba}^{(1)}_1}{\|\ba^{(1)}_1\|\,\|{\hat{\ba}}^{(1)}_1\|}~,
\end{equation}


\begin{theorem}\cite{2012arXiv1209.3215T}\label{theo_cribrnk2}
The Cram\'er-Rao induced bound (CRIB) on $\alpha_1^2$ for rank-2 tensor is given by
\begin{eqnarray}
\mbox{CRIB}(\ba_1)= \frac{\sigma^2}{\lambda_1^2 }\left(\frac{I_1-1}{1-h_1^2} +
\frac{(1-c_1^2)h_1^2}{1-h_1^2}\,\frac{y^2 + z - h_1^2z(z+1)}{(1 -
c_1y - h_1^2(z+1))^2+h_1^2(y + c_1z)^2}\right)\label{result31}
\end{eqnarray}
where $c_n = a^{(n) T}_1 a^{(n)}_2$, and 
\begin{eqnarray}
h_n &=& \prod_{2\leq k\ne n}^N c_k \qquad\mbox{for}\quad n=1,\ldots,N,\\
y &=& - c_1\,\sum_{n=2}^N \frac{h_n^2(1-c_n^2)}{c_n^2-h_n^2c_1^2},
\label{exp1}\\ z &=& \,\sum_{n=2}^N
\frac{1-c_n^2}{c_n^2-h_n^2c_1^2}~.\label{exp3}
\end{eqnarray}
\end{theorem}

\bibliographystyle{IEEEbib}
\bibliography{bibligraphy_thesis}

\end{document}